\documentclass[12pt]{amsart}
\usepackage{times}
\usepackage{amsfonts}
\usepackage{amsthm}
\usepackage{amsmath}
\advance \topmargin by -\headheight
\advance \topmargin by -\headsep
\evensidemargin \oddsidemargin
\newtheorem{theorem}{Theorem}[section]
\newtheorem{lemma}[theorem]{Lemma}
\newtheorem{corollary}[theorem]{Corollary}

\newtheorem{definition}[theorem]{Definition}

\begin{document}

\title{Fractal Entropies and Dimensions for Microstates Spaces}

\author{Kenley Jung}

\dedicatory{For Bill Arveson}

\address{Department of Mathematics, University of California,
Berkeley, CA 94720-3840,USA}

\email{factor@math.berkeley.edu}
\thanks{Research supported by the NSF Graduate Fellowship Program}

\begin{abstract} Using Voiculescu's notion of a matricial microstate we
introduce fractal dimensions and entropies for finite sets of selfadjoint
operators in a tracial von Neumann algebra.  We show that they possess
properties similar to their classical predecessors.  We relate the
new quantities to free entropy and free entropy dimension and show that a
modified version of free Hausdorff dimension is an algebraic invariant.  
We compute the free Hausdorff dimension in the cases where the set
generates a finite dimensional algebra or where the set consists of a
single selfadjoint.  We show that the Hausdorff dimension becomes additive
for such sets in the presence of freeness.  \end{abstract}

\maketitle

\section{Introduction} 

Voiculescu's microstate theory has settled some open questions in operator
algebras.  With it he shows in [10] that a free group factor does not have
a regular diffuse hyperfinite subalgebra (the first known kind with
separable predual).  Using similar techniques Ge shows in [2] that a free
group factor cannot be decomposed into a tensor product of two infinite
dimensional von Neumann algebras (again the first known example with
separable predual).  Both results rely upon the microstate theory
and the (modified) free entropy dimension. Free entropy dimension is a
number associated to an $n$-tuple of selfadjoint operators in a tracial
von Neumann algebra.  It is an analogue of Minkowski dimension and as such
it can be reformulated in terms of metric space packings.

   Unfortunately it is not known whether $\delta_0$ is an invariant of
von Neumann algebras, that is, if $\{b_1,\ldots,b_p\}$ is a set of
selfadjoint elements in $M$ which generates the same von Neumann algebra
as $\{a_1,\ldots,a_n\},$ then is it true that

\[ \delta_0(a_1,\ldots,a_n) = \delta_0(b_1,\ldots,b_p)? \]

\noindent The mystery of the invariance issue is this: how does the
asymptotic geometry of the microstate spaces reflect properties of the
generated von Neumann algebra of the $n$-tuple?

[5] shows that $\delta_0$ possesses a fractal geometric description in
terms of uniform packings.  Encouraged by this result we use microstates
to develop fractal geometric quantities for an $n$-tuple of selfadjoint
operators in a tracial von Neumann algebra.  For such an $n$-tuple
$z_1,\ldots, z_n$ we define several numerical measurements of their
microstate spaces, the most notable being the free Hausdorff dimension
of $z_1,\ldots,z_n.$ We denote this quantity by $\mathbb
H(z_1,\ldots,z_n).$ As in the classical case we have that $\mathbb
H(z_1,\ldots,z_n) \leq \delta_0(z_1,\ldots,z_n).$ For each $\alpha \in
\mathbb R_{+}$ we define an $\alpha$-free Hausdorff entropy for
$z_1,\ldots,z_n$ which we denote by $\mathbb
H^{\alpha}(z_1,\ldots,z_n).$ Hausdorff $n$-measure is a constant
multiple of Lebesgue measure and in our setting we have an analogous
statement: $\mathbb H^n(z_1,\ldots,z_n) = \chi(z_1,\ldots,z_n) +
\frac{n}{2}\log(\frac{2n}{\pi e}).$ A modified version of $\mathbb H$
denoted by $\overline {\mathbb H}$ turns out to be an algebraic
invariant.  We compute the free Hausdorff dimension of the $n$-tuple
when it generates a finite dimensional algebra or when $n=1.$ In both
cases the free Hausdorff and free entropy dimensions agree.  Using the
computations for a single selfadjoint, we show that if $\mathbb
H(z_1,\ldots, z_n) <1,$ then $\{z_1, \ldots, z_n\}^{\prime \prime}$ has
a minimal projection.  We view this as a microstates analogue of the
classical fact that a metric space with Hausdorff dimension strictly
less than $1$ must be totally disconnected.  Finally we show $\mathbb H$
becomes additive in the presence of freeness for the two aforementioned
$n$-tuples of random variables.

Our motivation in developing fractal dimensions for microstate spaces is
twofold.  Firstly, having other metric measurements for them may
eventually shed light on the connections between their asymptotic geometry
and the structure of the generated von Neumann algebras.  Secondly, it
provides another conceptual framework for the microstate theory.

Section 2 is a list of notation.  Section 3 reviews the definition of
classical Hausdorff dimension, then presents the free Hausdorff dimension
and entropy of an $n$-tuple and some of its basic properties.  The section
concludes with free packing and Minkowski-like entropies.  Section 4
introduces $\overline{\mathbb H},$ the modified version of $\mathbb H,$
and shows that $\overline {\mathbb H}$ is an algebraic invariant.  
Section 5 deals with the free Hausdorff dimension of finite dimensional
algebras.  Section 6 deals with the free Hausdorff dimension of single
selfadjoints.  Section 7 discusses various formulae for the free
Hausdorff dimension in the presence of freeness.  

\section{Preliminaries} Throughout this paper $M$ will be a von Neumann
algebra with a normal, tracial state $\varphi.$ $z_1, \ldots, z_n \in M$
are selfadjoint elements which generate a von Neumann algebra containing
the identity, and $\{s_i : i \in \mathbb N\}$ is always a semicircular
family in $M$ free with respect to $\{z_1,\ldots,z_n\}.$ We maintain the
notation for $\Gamma_R(:), \chi, \delta_0$ introduced in [9] and [10].  
$tr_k$ denotes the normalized trace on $M^{sa}_k(\mathbb C),$ the set of
$k \times k$ selfadjoint complex matrices, and $(M^{sa}_k(\mathbb C))^n$
is the set of $n$-tuples of elements in $M^{sa}_k(\mathbb C).$ $U_k$ is
the set of $k \times k$ complex unitaries.  $| \cdot|_2$ is the normalized
Hilbert-Schmidt norm on $M_k(\mathbb C)$ or $M$ induced by $tr_k$ or
$\varphi,$ respectively, or the norm on $(M_k(\mathbb C))^n$ given by
$|(x_1,\ldots,x_n)|_2 = \left (\sum_{i=1}^n tr_k(x_i^2) \right
)^{\frac{1}{2}}.$ Denote by vol Lebesgue measure on $(M^{sa}_k(\mathbb
C))^n$ with respect to the inner product on $(M^{sa}_k(\mathbb C))^n$
generated by the norm $\|(x_1,\ldots,x_n)\|_2 = (k \cdot \sum_{j=1}^n
tr_k(x_j^2))^{\frac{1}{2}}.$ For a metric space $(X,d)$ and $\epsilon >0$
write $P_{\epsilon} (X)$ for the maximum number of elements in a
collection of mutually disjoint open $\epsilon$ balls of $X.$ For a subset
$A$ of $X$ $|A|$ denotes the diameter of $A$ and $\mathcal
N_{\epsilon}(A)$ is the $\epsilon$ neighborhood of $A$ in $X.$

\section{Free Fractal Entropy and Dimensions}

Before defining a "free" Hausdorff dimension we recall classical 
Hausdorff dimension. 

\subsection{Hausdorff Dimension}

Suppose $(X,d)$ is a metric space.  For any $\epsilon, r >0$ define
$H^r_{\epsilon}(X)$ to be the infimum over all sums of the form $\sum_{j
\in J} |\theta_j|^r$ where $\langle \theta_j \rangle_{j \in J}$ is a
countable $\epsilon$-cover for X, i.e., $\langle \theta_j \rangle_{j \in
J}$ is countable collection of subsets of $X$ whose union contains $X$ and
for each $j$ $|\theta_j| \leq \epsilon.$ $H^r_{\epsilon}(\cdot)$ is an
outer measure on $X.$ Observe that if $\epsilon_1 < \epsilon_2,$ then
$H^r_{\epsilon_1}(X) \geq H^r_{\epsilon_2}(X).$ Hence, $H^r(A) =
\lim_{\epsilon \rightarrow 0} H^r_{\epsilon} (X) \in [0, \infty]$ exists.

If $r<s$ and $\epsilon >0,$ then for any countable $\epsilon$-cover
$\langle \theta_j \rangle_{j \in J}$ for $X,$

\[\epsilon^{s-r} \cdot \sum_{j \in J} |\theta_j|^r \geq \sum_{j \in J}
|\theta_j|^s.  \]

\noindent It follows that $H^r_{\epsilon}(X) \geq (\frac 
{1}{\epsilon})^{s-r}H_{\epsilon}^s(X).$  Taking a limit as $\epsilon 
\rightarrow 0$ shows that for any $K >0,$ $H^r(X) \geq K \cdot 
H^s(X).$  Consequently, there exists a nonnegative number $r$ for which 
$H^s(A) =0$ for all $s >r$ and $H^s(X) = \infty$ for all $s<r.$  
Formally, if $H^s(X) = \infty$ for all $s,$ then define $\dim_H(X) = 
\infty$.  Otherwise, define $\dim_H(X) = \inf \{r >0: H^r(X) = 0\}.$  
$\dim_H(X)$ is called the Hausdorff dimension of $X.$ 

	It turns out that $H^s(\cdot)$ generates a regular Borel measure
on $X$ called Hausdorff $s$-dimensional measure.  A subset $A$ of $X$ is
called an $s$-set if $0 < H^s(A) < \infty.$ It is clear that for every
$s$-set $A$ of $X$ $\dim_H(A) = s.$ However, it is possible for 
$A$ to have Hausdorff dimension $s$ and also satisfy $H^s(A) = 0$
or $H^s(A) = \infty.$

Manipulating the definitions one has that for any $S \in \mathbb R^d$
$\dim_H(S) \leq \dim_P(S)$ where $\dim_P(S)$ denotes the upper
Minkowski/uniform packing dimension of $S$ (see [1]).  There exist sets
$S$ for which the inequality is strict.

\subsection{Free Hausdorff Dimension}

In this subsection we define free Hausdorff 
dimension for $n$-tuples of selfadjoint elements in a von Neumann 
algebra and prove a few of its simple properties.

	In what follows, the Hausdorff and packing quantities on the
microstate spaces are taken with respect to the $| \cdot |_2$ metric
discussed in Section 2.  For any $m \in \mathbb N$ and $R, r, \epsilon,
\gamma >0$ define successively

\[\mathbb H_{\epsilon, R}^r(z_1, \ldots, z_n;m, \gamma) = \limsup_{k
\rightarrow \infty} \left (k^{-2} \cdot \log
[H_{\epsilon}^{rk^2}(\Gamma_R(z_1, \ldots,z_n;m,k,\gamma))] \right),\]

\[\mathbb H^r_{\epsilon,R}(z_1,\ldots,z_n) = \inf \{\mathbb
H^r_{\epsilon,R}(z_1,\ldots,z_n;m,\gamma): m\in \mathbb N, \gamma >0\}.\]

\noindent We now play the same limit games as in the classical case.  If 
$0 < \epsilon_1 < \epsilon_2,$ then for any $k,m,$ and $\gamma$ 

\[H^{rk^2}_{\epsilon_1}(\Gamma_R(z_1,\ldots,z_n;m,k,\gamma)) \geq 
H^{rk^2}_{\epsilon_2}(\Gamma_R(z_1,\ldots,z_n;m,k,\gamma)), \]

\noindent whence from the definitions, $\mathbb
H^r_{\epsilon_1,R}(z_1,\ldots,z_n)  \geq \mathbb
H^r_{\epsilon_2,R}(z_1,\ldots,z_n).$ Thus we define $\mathbb
H^r_R(z_1,\ldots,z_n)$ to be $\lim_{\epsilon \rightarrow 0} \mathbb
H^r_{\epsilon,R}(z_1,\ldots,z_n) \in [-\infty, \infty].$ If $r<s$ and
$1>\epsilon >0,$ then for any $k$

\[H_{\epsilon}^{rk^2}(\Gamma_R(z_1,\ldots,z_n;m,k,\gamma)) \geq
H_{\epsilon}^{sk^2}(\Gamma_R(z_1,\ldots,z_n;m,k,\gamma)) \cdot
\epsilon^{(r-s)k^2}. \]

\noindent Applying $k^{-2} \cdot \log$ to both sides, taking a $\limsup$
as $k \rightarrow \infty$ shows that for any $m \in \mathbb N$ and 
$R, \gamma
> 0,$

\[\mathbb H_{\epsilon,R}^{r}(z_1,\ldots,z_n;m,\gamma) \geq \mathbb
H_{\epsilon,R}^{s}(z_1,\ldots,z_n;m,\gamma) +(s-r)| \log \epsilon|. \]

\noindent Taking infimums over $m$ and $\gamma$ followed by a limit as
$\epsilon \rightarrow 0$ we have for any $R, K >0$

\[ \mathbb H^r_R(z_1,\ldots,z_n) \geq \mathbb H^s_R(z_1,\ldots,z_n) + K. \]  

\begin{definition} The free 
Hausdorff $r$-entropy of $z_1,\ldots,z_n$ is 

\[\mathbb H^r(z_1,\ldots, z_n) = \sup_{R>0} \mathbb H^r_R(z_1,\ldots,z_n).\]
\end{definition}
 
\begin{lemma} If $r > \delta_0(z_1,\ldots,z_n),$ then $\mathbb
H^r(z_1,\ldots,z_n) = - \infty.$ \end{lemma} \begin{proof} For any $R,
\epsilon, \gamma>0$ and $m,k \in \mathbb N$ it is clear that \[\log \left
[ H^{rk^2}_{4\epsilon}(\Gamma_R(z_1,\ldots,z_n;m,k,\gamma))  \right] \leq
\log(P_{\epsilon}(\Gamma_R(z_1,\ldots,z_n;m,k,\gamma))) + rk^2 \cdot \log
(4 \epsilon).\]

\noindent Consequently by [5],  

\begin{eqnarray*} \mathbb H^r_R(z_1,\ldots,z_n) = \lim_{\epsilon 
\rightarrow 0} \mathbb
H^r_{4\epsilon,R}(z_1,\ldots,z_n)  & = &  \limsup_{\epsilon \rightarrow 0}
\mathbb H^r_{4\epsilon,R}(z_1,\ldots, z_n) \\ & \leq & \limsup_{\epsilon 
\rightarrow 0} \mathbb P_{
\epsilon}(z_1,\ldots,z_n) + r \cdot \log 4 \epsilon \\ & = &
\limsup_{\epsilon \rightarrow 0} \left( \frac {\mathbb P_{ \epsilon} (z_1,   
\ldots, z_n)}{| \log \epsilon|} \cdot |\log \epsilon| + r \cdot \log
4 \epsilon \right) \\ & \leq & \limsup_{\epsilon \rightarrow 0} 
\left ( r \cdot | \log  \epsilon| + r \cdot \log 4 \epsilon \right)
\\
& = & r \cdot \log 4.
\end{eqnarray*}

\noindent Hence, $\mathbb H^r(z_1,\ldots,z_n) < \infty$ for all
$r > \delta_0(z_1,\ldots,z_n)$ and the result follows.
\end{proof}

\begin{definition} If $\{z_1,\ldots,z_n\}$ has finite dimensional
approximants, then $\mathbb H(z_1,\ldots,z_n) = \inf \{ r>0:  \mathbb
H^r(z_1,\ldots,z_n) = -\infty\}.$ Otherwise, define $\mathbb
H(z_1,\ldots,z_n) = -\infty.$ $\mathbb H(z_1,\ldots,z_n)$ is called the
free Hausdorff dimension of $z_1,\ldots,z_n.$ \end{definition}

\begin{definition} For $s \geq 0$ $\{z_1,\ldots,z_n\}$ is an $s$-set if 
$-\infty < \mathbb H^s(z_1,\ldots,z_n) < \infty.$
\end{definition}

By definition if $\{z_1,\ldots,z_n\}$ is an $s$-set, then $\mathbb 
H(z_1,\ldots,z_n) = s.$  

Classical uniform packing dimension dominates Hausdorff dimension and it 
is not surprising that in our setting we have the analogous statement by 
Lemma 3.2 and Definition 3.3:

\begin{corollary} $\mathbb H(z_1,\ldots,z_n) \leq 
\delta_0(z_1,\ldots,z_n).$
\end{corollary}

\begin{lemma} If $y_1,\ldots,y_p$ are self-adjoint elements in 
$\{z_1,\ldots,z_n\}^{\prime \prime},$ then for any $r>0$

\[ \mathbb H^r(z_1,\ldots,z_n) \leq \mathbb 
H^r(z_1,\ldots,z_n,y_1,\ldots,y_p).\]  
\end{lemma}

\begin{proof} Withough loss of generality assume that the $z_i$ have 
finite dimensional approximants.  Suppose $R$ exceeds the operator norms 
of the $z_i$ or $y_j.$  Given $m \in \mathbb N$ and $\epsilon, \gamma >0$ 
there exist by Lemma 4.1 of [4] $m_1 \in \mathbb N,$ $\gamma_1>0$ and 
polynomials $f_1,\ldots,f_p$ in $n$ noncommutative variables such that if 
$(x_1,\ldots, x_n) \in \Gamma_R(z_1,\ldots,z_n;m_1,k,\gamma_1)$ then

\[ (x_1,\ldots,x_n,f_1(x_1,\ldots,x_n),\ldots,f_p(x_1,\ldots,x_n)) \in 
\Gamma_R(z_1,\ldots,z_n,y_1,...,y_p;m,k,\gamma).\]

\noindent For each $k$ this map from
$\Gamma_R(z_1,\ldots,z_n;m_1,k,\gamma_1)$ to
$\Gamma_R(z_1,\ldots,z_n,y_1,\ldots,y_p;m,k,\gamma)$ defined by sending
$(x_1,\ldots,x_n)$ to $(x_1,\ldots,x_n, f_1(x_1,\ldots,x_n),\ldots,
f_p(x_1,\ldots,x_n))$ increases distances with respect to $|\cdot|_2.$
Hence

\[ \mathbb H^r_{\epsilon,R}(z_1,\ldots,z_n;m_1,\gamma_1) \leq \mathbb 
H^r_{\epsilon,R}(z_1,\ldots,z_n,y_1,\ldots,y_p;m,\gamma).\]

\noindent This being true for any $m,\gamma, \epsilon,$ and $R$ as
stipulated, the results follows.  \end{proof}

In [10] it was shown that $\chi(z_1,\ldots,z_n) > - \infty \Rightarrow
\delta_0(z_1,\ldots,z_n) =n.$ This is the noncommutative analogue of the
fact that if a Borel set $S \subset \mathbb R^d$ has nonzero Lebesgue
measure, then its uniform packing dimension is $d.$ One can replace
"uniform packing" in the preceding sentence with "Hausdorff," and we
record its analogue, after making a simple remark about free Hausdorff
entropy.

   	 The following equation says that free entropy is a normalization
of free Hausdorff $n$-entropy and echoes asymptotically in statement and
proof the fact that Lebesgue measure is a scalar multiple of Hausdorff
dimension.

\begin{lemma} $ \mathbb H^n(z_1,\ldots,z_n) = \chi(z_1,\ldots,z_n) + K$ 
where $K =
\frac{n}{2} 
\log \left( \frac{2n}{\pi e} \right).$
\end{lemma}

\begin{proof}We can clearly assume that $\{z_1,\ldots,z_n\}$ has finite
dimensional approximants.  First we show that the left hand side of the
equation is greater than or equal to the right hand side.  Suppose that
$\epsilon, \gamma>0, m,k \in \mathbb N, R > \max_{1 \leq j \leq
n}\{\|z_j\|\}.$ Suppose $\langle \theta_j \rangle_{j \in J}$ is a cover of
$\Gamma_R(z_1,\ldots,z_n;m,k,\gamma).$ Because any set is contained in a
closed convex set of equal diameter we may assume that the $\theta_j$ are
closed and convex.  In particular they are Borel sets and the
isodiametric inequality yields \[ \sum_{j \in J} |\theta_j|^{nk^2} \geq
\frac{2^{nk^2} \Gamma(\frac{nk^2}{2}+1)}{\sqrt{\pi k}^{nk^2}} \cdot
\sum_{j \in J} \text{vol}(\theta_j) \geq
\frac{2^{nk^2}\Gamma(\frac{nk^2}{2}+1)}{\sqrt{\pi k }^{nk^2}} \cdot 
\text{vol}
(\Gamma_R(z_1,\ldots,z_n;m,k,\gamma)). \]

\noindent Thus, 
$\mathbb H^{n}_{\epsilon}(\Gamma_R(z_1,\ldots,z_n;m,\gamma))$
dominates
\begin{eqnarray*} &  &   \limsup_{k\rightarrow \infty} \left( k^{-2} \cdot
\log[\text{vol}(\Gamma_R(z_1,\ldots,z_n;m,k,\gamma))] + k^{-2} \cdot
\log\left(\Gamma\left(\frac{nk^2}{2}+1\right)\right) - 
\frac{n}{2} 
\log \frac { \pi k }{4} \right) \\  & \geq &  \limsup_{k \rightarrow 
\infty} \left( k^{-2} \cdot \log
[\text{vol}(\Gamma_R(z_1,\ldots,z_n;m,k,\gamma))] + \frac{n}{2} \log 
\left(
\frac{nk^2}{2e} \right) - \frac{n}{2} \log \frac {\pi k}{4} \right)   
\\ & \geq & \chi_{R}(z_1,\ldots,z_n) + K.
\end{eqnarray*}

\noindent The above is a uniform lower bound for any $R, m, \gamma,$ and 
$\epsilon$ so 
\[\mathbb H^n(z_1,\ldots,z_n) \geq \chi(z_1,\ldots,z_n) +K.\]

For the reverse inequality again assume $\epsilon, \gamma, m,k,$ and $R$ 
are as before.  Given $k$ use Vitali's covering lemma to find an 
$\epsilon$-cover $\langle \theta_j \rangle_{j \in J}$ for 
$\Gamma_R(z_1,\ldots, z_n;m,k,\gamma)$ such that each $\theta_j$ is a 
closed ball and 

\[ \sum_{j \in J} \text{vol}(\theta_j) < 2 \cdot
\text{vol}(\Gamma_R(z_1,\ldots,z_n;m,k,\gamma)).\]

\noindent By definition 

\[ \frac{(\pi k)^{\frac {nk^2}{2}}}{2^{nk^2}\Gamma(\frac{nk^2}{2} +1)} 
\cdot 
H^{nk^2}_{\epsilon}(\Gamma_R(z_1,\ldots,z_n;m,k, \gamma)) \leq 
2 \cdot \text{vol}(\Gamma_R(z_1,\ldots,z_n;m,k,\gamma)).\]

\noindent Thus $\mathbb H^n_{\epsilon, R}(z_1,\ldots,z_n;m,\gamma)$ is 
dominated by 

\begin{eqnarray*} & & \limsup_{k \rightarrow \infty} \left (
k^{-2} \cdot \log (\text{vol} (\Gamma_R(z_1,\ldots,z_n;m,k,\gamma))) + 
\frac{n}{2} \log \left (\frac {4}{\pi k} \right) + k^{-2} \cdot \log 
\left(\Gamma(\frac{nk^2}{2}+1)\right) \right) \\ & = & 
\limsup_{k \rightarrow \infty} \left(k^{-2} \cdot 
\log (\text{vol}(\Gamma_R(z_1,\ldots,z_n;m,k,\gamma))) + \frac {n}{2} 
\cdot \log 
\left (\frac 
{nk^2}{2e} \right) + \frac{n}{2} \log \left( \frac{4}{\pi k} \right) 
\right) \\ & 
= & \chi_{R} (z_1,\ldots,z_n;m,\gamma) + K.
\end{eqnarray*}
\noindent $\mathbb H^n_{\epsilon, R}(z_1,\ldots,z_n) \leq 
\chi_R(z_1,\ldots,z_n) + K.$  It follows that 

\[\mathbb H^n(z_1,\ldots,z_n) \leq \chi(z_1,\ldots,z_n) +K.\]           
\end{proof}

We now have:

\begin{corollary} If $\chi(z_1,\ldots,z_n) > - \infty,$ then 
$\{z_1, \ldots, z_n\}$ is an $n$-set.  In particular $\mathbb 
H(z_1,\ldots,z_n) = n.$
\end{corollary}

	Observe that if $s_1, \ldots, s_n$ is a free semicircular family, 
then $\chi(s_1,\ldots,s_n) > - \infty$ by [9] whence $\mathbb 
H(s_1,\ldots,s_n) = \delta_0(s_1,\ldots,s_n) = n$ by Corollary 3.8.

It may have crossed the reader's mind why we did not prove or in
the very least make a remark about a subadditive property for $\mathbb
H.$ At this point we recall a difference between classical
Hausdorff dimension and Minkowski dimension when taking Cartesian
products.  Suppose $S \subset \mathbb R^m$ and $T \subset \mathbb R^n$
are Borel sets, $\dim_H(\cdot)$ is Hausdorff dimension, and
$\dim_M(\cdot)$ denotes Minkowski dimension.  It is easy to see that

\[ \dim_M(S \times T) \leq \dim_M(S) + \dim_M(T).\]

\noindent On the other hand with some work (see [1]) one shows

\[ \dim_H(S \times T) \geq \dim_H(S) + \dim_H(T). \]

\noindent Strictness of the above inequality can occur.  In fact, there
exist sets $S, T \subset \mathbb R$ with Hausdorff dimension 0 such that
$\dim_H (S \times T) =1$ (see [1]).  We do not know if there exist 
sets of self-adjoints satisfying a similar inequality.  

In general we have $\dim_H(S \times T) \leq \dim_H(S) + \dim_M(T).$ The
classical proof does not immediately provide a proof for the microstates
situation.  The obstruction occurs when one fixes the parameter $\epsilon$
and finds that the convergence of the $\epsilon$ packing number of the
$k\times k$ matricial microstates as $k \rightarrow \infty$ may depend too
heavily upon the choice of $\epsilon$ and thus grow too slow for good
control as one pushes $\epsilon$ to $0.$ In some cases, however, one can
obtain strong upper bounds where for small enough $\epsilon$ the rate of
convergence of the $k$-dimensional quantities behaves appropriately so
that the inequality $\mathbb H(y_1,\ldots,y_m, z_1,\ldots,z_n) \leq
\mathbb H(y_1,\ldots,y_m) + \delta_0(z_1,\ldots,z_n)$ holds.  In
particular, the inequality will occur when $\{z_1,\ldots,z_n\}$ generates
a hyperfinite von Neumann algebra or when it can be partitioned into a
free family of sets each of which generates a hyperfinite von Neumann
algebras.

\subsection{Free Entropies for $\delta_0$}

Although Hausdorff dimension and measure can provide metric information
about sets they are often difficult to compute (in particular, finding
sharp lower bounds is a problem).  A machine which sometimes allows for
easier computations is Minkowski content.  It is a numerical measurement
related to Minkowski dimension in the same way that Hausdorff measure is
related to Hausdorff dimension.

Suppose $X \subset \mathbb R^d.$  For $r >0$ define $M^r(X)$ by

\[ \limsup_{\epsilon \rightarrow 0} \frac {\lambda(\mathcal 
N_{\epsilon}(X))}{\lambda (B_{\epsilon}^{d-s})} \]

\noindent where $\lambda$ is Lebesgue measure on $\mathbb R^d$ and 
$B_{\epsilon}^{d-s}$ is the ball of radius $\epsilon$ centered at the 
origin in $\mathbb R^{d-s}.$  As with 
$M^r(\cdot)$ we have that $M^r(X) \geq K \cdot M^s(X)$ for any $r < s$ and 
$K>0.$  Hence there exists a nonnegative number $r$ for which $M^s(X) = 
0$ for $r<s$ and $M^s(X) = \infty$ for $r>s.$  

This number $r$ turns out to be the Minkowski dimension of $X.$  $M^r(X)$ 
is called the Minkowski content of $X$ and 
provides a measurement of the size of $X.$  We can also define a packing 
quantity related to $M^r(X),$ $P^r(X),$ by

\[ \limsup_{\epsilon \rightarrow 0} P_{\epsilon}(X) \cdot 
(2\epsilon)^r.\]  

\noindent As before there exists a unique $r \geq 0$ for which $P^s(X) = 
0$ if $r <s$ and $P^s(X) = \infty$ if $r >s.$  Again this unique $r$ turns 
out to be the Minkowski dimension of $X.$

Unlike the Hausdorff construction neither $M^r$ nor $P^r$ are measures.  
In fact, they have the unpleasant property (from the classical viewpoint) 
that their values of a set and its closure are the same.  Hence, 
$M^1(\mathbb Q) = P^1(\mathbb Q) = \infty.$  On the other hand  
$H^1(\mathbb Q)=0$ (although $H^0(\mathbb Q) = \infty$).

In the following $\mathcal N_{\epsilon}$ and will be taken 
with respect to the $|\cdot|_2$ metrics.     

\begin{definition} For any $m \in \mathbb N$ and $R, \alpha, \gamma, 
\epsilon >0$ define successively,

\begin{eqnarray*} \mathbb M^{\alpha}_{\epsilon, R}(z_1,\ldots,z_n;m,
\gamma) & = & \limsup_{k \rightarrow \infty} k^{-2} \cdot \log (
\text{vol} (\mathcal N_{\epsilon}(\Gamma_R(z_1,\ldots,z_n;m,k,\gamma)))) +
\frac {n}{2} \log k \\ & & + |\log (\epsilon^{n- \alpha})| \end{eqnarray*}

\[ \mathbb M^{\alpha}_{\epsilon, R}(z_1,\ldots,z_n) = \inf 
\{\mathbb M^{\alpha}_{\epsilon,R}(z_1,\ldots,z_n;m,\gamma): m\in \mathbb 
N, \gamma >0 \}, \]

\[ \mathbb M^{\alpha}_{R}(z_1,\ldots,z_n) = \limsup_{\epsilon 
\rightarrow 0} \mathbb 
M^{\alpha}_{\epsilon,R}(z_1,\ldots,z_n) \]

\[ \mathbb M^{\alpha}(z_1,\ldots,z_n) = \sup_{R> 0}
\mathbb M^{\alpha}_{R}(z_1,\ldots,z_n).\]

\noindent We call $\mathbb M^{\alpha}(z_1,\ldots,z_n)$ the free Minkowski 
$\alpha$-entropy of $\{z_1,\ldots,z_n\}.$  

\end{definition}

Recalling the definition of $\mathbb P_{\epsilon}(z_1,\ldots,z_n)$ in [5] 
we also make the following:

\begin{definition} For $\alpha >0$ the free packing $\alpha$-entropy of 
$\{z_1,\ldots,z_n\}$ is 

\[ \mathbb P^{\alpha}(z_1,\ldots,z_n) = \sup_{R >0} \mathbb 
P^{\alpha}_R(z_1,\ldots,z_n) \] 

\noindent where

\[ \mathbb P^{\alpha}_R(z_1,\ldots,z_n) = \limsup_{\epsilon \rightarrow 0}
\mathbb P_{\epsilon,R}(z_1,\ldots,z_n) + \alpha \cdot \log 2 \epsilon . \]
\end{definition}

The following is a easy and we omit the proof.  It relates $\mathbb
H^{\alpha}, \mathbb P^{\alpha},$ the free entropy of an
$\epsilon$-semicircular perturbation, and $\mathbb M^{\alpha}$ and shows
that the latter three give the same information modulo universal
constants.

\begin{lemma} For any $\alpha >0$ 

\begin{eqnarray*} \mathbb H^{\alpha}(z_1,\ldots,z_n) - \alpha \log 2 & 
\leq & \mathbb P^{\alpha}(z_1,\ldots, z_n) \\ & \leq &
\limsup_{\epsilon \rightarrow 0} \left [ \chi(z_1 + \epsilon s_1, \ldots,
z_n + \epsilon s_n:s_1, \ldots, s_n) + (n - \alpha ) |\log \epsilon|
\right] \\ & & + \alpha \cdot \log (4\sqrt {n}) - \chi(s_1,\ldots,s_n) \\ 
& \leq
& \mathbb M^{\alpha}(z_1,\ldots,z_n) + (n-\alpha) \log \sqrt{n} + \alpha 
\cdot \log (4\sqrt{n}) -
\frac{n}{2} \log (2 \pi e)  \\ & \leq & \mathbb
P^{\alpha}(z_1,\ldots,z_n) + \alpha \log 4\sqrt{n} + n \log 4. 
\end{eqnarray*}
\end{lemma}

\section{Modified Free Hausdorff Dimension and Algebraic Invariance}

Thus far we cannot prove that $\mathbb H$ is an algebraic invariant and
towards this end we introduce a technical modification of $\mathbb H.$ For
a metric space $(X,d),$ and $0<\delta < \epsilon$ a $(\delta<
\epsilon)$-cover $\langle \theta_j \rangle_{j \in J}$ is a countable cover
of $X$ such that for each $j$ $\delta \leq |\theta_j| \leq \epsilon.$ For
$r >0$ define $H^r_{(\delta < \epsilon)}(X)$ to be the infimum over all
sums of the form $\sum_{j \in J} |\theta_j|^s$ where $\langle \theta_j
\rangle_{j \in J}$ is a $(\delta < \epsilon)$-cover of $X.$ As before, for
$\delta < \epsilon_1 < \epsilon_2$ and $s > r \geq 0$ $H^r_{(\delta <
\epsilon _1)}(X) \geq H^r_{(\delta < \epsilon_2)}(X)$ and $H^r_{(\delta <
\epsilon)}(X) \geq H^s_{(\delta < \epsilon)} \cdot \frac
{1}{\epsilon^{s-r}}.$

Taking all Hausdorff quantities with respect to the $| \cdot |_2$ metric,
define successively for any $m \in \mathbb N,$ and
$L,R,r,\epsilon,\gamma>0$ with $L \sqrt{\gamma} < \epsilon$

\[ \overline{\mathbb H}^r_{L, \epsilon,R}(z_1,\ldots,z_n;m,\gamma) =
\limsup_{k \rightarrow \infty} \left[k^{-2} \cdot \log 
H^{rk^2}_{(L\sqrt{\gamma} <
\epsilon), R}(z_1,\ldots,z_n;m,\gamma)\right],\]

\[ \overline{\mathbb H}^r_{L, \epsilon, R}(z_1,\ldots z_n) = \inf \{
\overline{\mathbb H}^r_{L , \epsilon, R}(z_1,\ldots,z_n;m,\gamma) : m\in
\mathbb N, \gamma >0 \}, \]

\[ \overline{\mathbb H}^r_{\epsilon,R}(z_1,\ldots, z_n) = \sup_{L>0}
\overline{\mathbb H}^r_{L,\epsilon, R}(z_1,\ldots, z_n).\]

\noindent As before $\overline {\mathbb
H}^r_{\epsilon_1,R}(z_1,\ldots,z_n) \geq \overline {\mathbb
H}^r_{\epsilon_2,R}(z_1,\ldots,z_n)$ for $\epsilon_1 < \epsilon_2$ so that
$\lim_{\epsilon \rightarrow 0} \overline{\mathbb
H}^r_{\epsilon,R}(z_1,\ldots,z_n) \in [-\infty, \infty]$ exists.  Write
$\overline{\mathbb H}^r_R(z_1,\ldots,z_n)$ for this limit and
$\overline{\mathbb H}^r(z_1,\ldots,z_n) = \sup_{R>0} \overline{\mathbb
H}^r_R(z_1,\ldots,z_n).$ Obviously $\mathbb H^r(z_1,\ldots,z_n) \leq
\overline{\mathbb H}^r(z_1,\ldots,z_n).$ For any $r < s$ and $K>0$ we
have $\overline{\mathbb H}^{r}(z_1,\ldots, z_n) \geq \overline{\mathbb
H}^{s}(z_1,\ldots,z_n) +K.$

	We have the analogous result of Lemma 3.2 provided we know that
uniform packings by open $\epsilon$-balls of microstate spaces generate
$(L\sqrt{\gamma}< \epsilon)$-covers for $\gamma$ sufficiently small.  
This is not immediate for a priori an $\epsilon$ ball in a microstate
space could have diameter much smaller than $\epsilon,$ possibly even $0$
(in this case the ball consists of just a single point).  But
path-connectedness of $U_k$ and a simple point set topology argument
allows us to deduce that for large dimensions the microstate spaces are
rich enough so that the diameter of any $\epsilon$ ball with microstate
center is at least $\epsilon$:

\begin{lemma} Suppose $\{z_1,\ldots,z_n\}$ generates a von Neumann algebra
not equal to $\mathbb C I$ and $R > \max\{\|z_i\|\}_{1\leq i \leq n}.$
There exist $\epsilon_0, \gamma >0,$ and $m, N \in \mathbb N$ such that if
$\epsilon_0 > \epsilon >0, k \geq N,$ and $(x_1,\ldots,x_n) \in
\Gamma_R(z_1,\ldots,z_n;m,k,\gamma),$ then there is a $Y \in
\Gamma_R(z_1,\ldots,z_n;m,k,\gamma)$ with

\[|Y -(x_1,\ldots,x_n)|_2 = \epsilon.\] \end{lemma} 

\begin{proof} There exists some $i$ such that $z_i \notin \mathbb C I.$
Without loss of generality we can assume that $z= z_1 \notin \mathbb C I.$
Hence by [10] $\delta_0(z) >0.$ Set $\beta = \frac{\delta_0(z)}{2}.$ By
[5] find some $\epsilon_0$ satisfying $1/40 > \epsilon_0 >0$ and

\[ \mathbb P_{20\epsilon_0,R}(z) > \beta \cdot |\log 20 \epsilon_0|.\]

\noindent Thus by regularity of $\chi$ for a single self-adjoint and 
[5] there exist $m, N \in \mathbb N$ and $\gamma >0$ such that 
if $k >N,$ then 

\[ k^{-2} \cdot \log (P_{20 \epsilon_0}(\Gamma_R(z;m,k,\gamma))) > \beta
\cdot |\log 20 \epsilon_0|.\]

\noindent By [4] we may choose $m$ and $\gamma$ so that if $x, y \in 
\Gamma_R(z;m,k,\gamma),$ then there exists a $u \in U_k$ such that $|uxu^* 
- y |_2 < \epsilon_0.$  

Now suppose that $(x_1,\ldots,x_n) \in 
\Gamma_R(z_1,\ldots,z_n;m,k,\gamma)$ with $k > \max \{N, 1/ 
\sqrt{\beta}\}.$  Suppose also that 
$0<\epsilon<\epsilon_0$ and $B_{\epsilon}$ is the open ball of $|\cdot|_2$ 
-radius $\epsilon$ centered at $(x_1,\ldots,x_n).$  Assume by 
contradiction that $\partial B_{\epsilon} \bigcap 
\Gamma_R(z_1,\ldots,z_n;m,k,\gamma) = \emptyset.$  Write 
$U(x_1,\ldots,x_n)$ for the set of all $n$-tuples of the form 
$(ux_1u^*,\ldots,ux_nu^*)$ for $u \in U_k.$  Clearly $\partial 
B_{\epsilon} \bigcap U(x_1,\ldots,x_n) = \emptyset.$  $U(x_1,\ldots,x_n)$ 
is a path-connected set so this implies that $U(x_1,\ldots,x_n)$ is 
contained in $B_{\epsilon}.$  Consequently, $U(x_1) = \{ ux_1u^*: u \in 
U_k\}$ is contained in the open ball of $|\cdot |_2$-radius $\epsilon$ 
with center $x_1.$  

On the other hand the selection of $m$ and $\gamma$ imply 
$\Gamma_R(z;m,k,\gamma) \subset \mathcal N_{\epsilon_0}(U(x_1)).$  
Combined with the estimate of the first paragraph we
have

\[ P_{4\epsilon_0}(U(x_1)) \geq P_{20 \epsilon_0}(\mathcal 
N_{\epsilon_0}(U(x_1))) \geq \left(\frac{1}{20 \epsilon_0}\right)^{\beta 
k^2} >2.\]

\noindent Thus one can find two points in $U(x_1)$ whose $|\cdot|_2$
distance from one another is greater than or equal to $4\epsilon_0
>4\epsilon.$ It follows that $U(x_1)$ cannot possibly be covered by the
open ball of $|\cdot|_2$-radius $\epsilon.$ This is absurd.  There must
exist some $Y \in \Gamma_R(z_1,\ldots,z_n;m,k,\gamma)$ with $|Y -
(x_1,\ldots,x_n)|_2 = \epsilon.$ \end{proof}
  
Lemma 4.1 with the proof of Lemma 3.2 show that if $r >
\delta_0(z_1,\ldots,z_n),$ then $\overline{\mathbb H}^{r}(z_1,\ldots,z_n)
= - \infty,$ provided the $z_i$ generate a nontrivial von Neumann algebra.
Otherwise they generate $\mathbb C I$ and then it's clear that for $r >0$ 
$\overline{\mathbb H}^r(z_1,\ldots,z_n) = -\infty.$

\begin{definition} The modified free Hausdorff dimension of
$\{z_1,\ldots,z_n\}$ is

\[\overline{\mathbb H}(z_1,\ldots,z_n) = \inf \{r >0: \overline{\mathbb
H}^r(z_1,\ldots,z_n) = -\infty\}.\]
\end{definition}

Immediately we observe that:

\begin{corollary} $\mathbb H(z_1,\ldots,z_n) \leq \overline{\mathbb
H}(z_1,\ldots,z_n) \leq \delta_0(z_1, \ldots,z_n).$
\end{corollary}

We now come to the primary result concerning $\overline{\mathbb H}.$
Viewing polynomials as Lipschitz maps when restricted to norm bounded
sets, the following is simply a corollary of the fundamental fact that
fractal dimensions are preserved under bi-Lipschitz maps.

\begin{lemma} If $\{y_1,\ldots,y_p\}$ and $\{z_1,\ldots,z_n\}$ are
sets of selfadjoint elements in $M$ which generate the same algebra,
then

\[ \overline{\mathbb H}(y_1,\ldots,y_p) = \overline{\mathbb
H}(z_1,\ldots,z_n).\]

\end{lemma}

\begin{proof} Set $Y = \{y_1,\ldots,y_p\}, Z = \{z_1,\ldots,z_n\}.$
Write $\Gamma_R(Y;m,k,\gamma)$ for $\Gamma_R(y_1,\ldots,y_p,m,k,\gamma).$
We interpret quantities like $\overline{\mathbb H}(Y)$ in the obvious way.  
Similarly for $Z.$ It suffices to show that $\overline {\mathbb H}(Y)  
\leq \overline {\mathbb H}(Z).$ Find $n$ polynomials $f_1,\ldots,f_n$ in
$p$ noncommuting variables such that $f_j(y_1,\ldots,y_p) = z_j$ for each
$j.$ Similarly find $p$ polynomials $g_1, \ldots, g_p$ in $n$ noncommuting
variables such that for each $i$ $g_i(z_1,\ldots,z_n)=y_i.$ For any
$(a_1,\ldots,a_p) \in (M^{sa}_k(\mathbb C))^p$ define

\[\Phi(a_1,\ldots,a_p) = (f_1(a_1,\ldots,a_p), \ldots
f_n(a_1,\ldots,a_p)).\]

\noindent Similarly for any $(b_1,\ldots,b_n) \in (M^{sa}_k(\mathbb C))^n$
define \[\Psi(b_1,\ldots,b_n) = (g_1(b_1,\ldots,b_n),
\ldots,g_p(b_1,\ldots,b_n)).\]

\noindent We can arrange it so that $\Phi$ takes $p$-tuples of 
$k \times k$ selfadjoints to $n$-tuples of $k \times k$ self-adjoints for 
any $k.$  Similarly for $\Psi.$

Suppose $R >0, 1 > \epsilon >0.$ There exists a constant $L>1$ dependent
only upon $\Phi$ and $R$ such that for any $x, y \in ((M^{sa}_k(\mathbb
C))_R)^p,$ $\Phi(x) \in ((M^{sa}_k(\mathbb C))_{LR})^n$ and $|\Phi(x) -
\Phi(y)|_2 \leq L |x-y|_2.$ Also, there exist $K>0$ and $m_0 \in \mathbb
N$ dependent only upon the polynomial $\Psi \circ \Phi$ and $R$ such that
for any $\gamma >0$ and $x \in \Gamma_R(Y;m_0,k,\gamma)$

\[ |(\Psi\circ \Phi)(x) - x|_2 < K \sqrt{\gamma}.\] 

\noindent Now suppose $m \in \mathbb N$ and $\gamma >0$ with
$2KL\sqrt{\gamma} < \epsilon.$ Choose $m_1 \in \mathbb N$ and $\gamma 
>\gamma_1 >0$ such that for any $k \in \mathbb N$

\[ \Phi(\Gamma_R(Y;m_1,k,\gamma_1)) \subset \Gamma_{LR}(Z;m,k,\gamma) 
\]

\noindent There exists a constant $L_1 >0$ dependent upon $LR$ and $\Psi$
such that for any $a,b \in ((M^{sa}_k(\mathbb C))_{LR})^n$

\[ |\Psi(a) - \Psi(b)|_2 \leq L_1 \cdot |a-b|_2.\]

Suppose $\langle \Theta_s \rangle_{s \in S}$ is a $(2KL\sqrt{\gamma} <
\epsilon)$-cover
for $\Gamma_{LR}(Z;m,k,\gamma).$ Define $S_1$ to be the subset of $S$
consisting of those elements $i$ such that $\Theta_i$ has nontrivial
intersection with $\Phi(\Gamma_R(Y;m_1,k,\gamma_1)).$ For each $i \in S_1$
define $B_i$ to be an open ball of $|\cdot|_2$-radius $|\Theta_i|$ with
center in $\Theta_i \bigcap \Phi(\Gamma_R(Y;m_1,k,\gamma_1)).$ Consider
for each such $i$ the center of the ball $B_i.$ Take one element $x$ in
the preimage of this center under $\Phi$ as a map restricted to
$\Gamma_R(Y;m_1,k,\gamma_1).$ In the metric space
$(\Gamma_R(Y;m_1,k,\gamma_1), | \cdot |_2)$ consider $B(x, |\Theta_i|/L),$
the open ball of radius $|\Theta_i|/L$ with center $x.$ $B(x, |\Theta_i|/L
)$ has diameter no less that $|\Theta_i|/L$ (by Lemma 4.1 with appropriate
restrictions on $\epsilon, m, k, \gamma$ and $R$) and clearly it lies in
$\Phi^{-1}(\Theta_i)$ so that $|\Phi^{-1}(\Theta_i)| \geq |\Theta_i|/L
\geq 2K \sqrt{\gamma}.$ For any $y,w \in \Gamma_R(Y;m_1,k,\gamma_1))$ we 
have

\[  L_1 \cdot |\Phi(y) - \Phi(w)|_2 \geq 
|(\Psi \circ \Phi)(y) - (\Psi \circ \Phi)(w)|_2 \geq |y - w|_2 - 2 K 
\sqrt{\gamma_1}\]

\noindent whence it follows that $L_1 \cdot |\Theta_i| > 
|\Phi^{-1}(\Theta_i)| 
-2 K \sqrt{\gamma_1} \Rightarrow (L_1+1) \epsilon > \Phi^{-1}(\Theta_i).$

	$\langle |\Phi^{-1}(\Theta_i)| \rangle_{i \in S_1}$ is a $(2K
\sqrt{\gamma}< (L_1+1) \epsilon)$ cover for $\Gamma_R(Y;m_1,k,\gamma_1).$
Estimate:

\begin{eqnarray*} H_{\left(2K \sqrt{\gamma_1}< 
(L_1+1)\epsilon \right)}^{rk^2}(\Gamma_R(Y;m_1,k,\gamma_1)) \leq \sum_{i 
\in S_1} |\Phi^{-1}(\Theta_i)|^{rk^2} & \leq & \sum_{i \in S_1} 
(L_1|\Theta_i| + 2K \sqrt{\gamma_1})^{rk^2} \\ & \leq & \sum_{i \in S_1} 
((L_1+1)|\Theta_i|)^{rk^2} \\ & \leq & (L_1+1)^{rk^2} \cdot \sum_{i \in S} 
|\Theta_i|^{rk^2}.\\ \end{eqnarray*}

\noindent This being true for any $(2KL\sqrt{\gamma} < \epsilon)$-cover 
$\langle \Theta_i \rangle_{s \in S}$ for $\Gamma_{LR}(Z;m,k,\gamma)$ it
follows that

\[H^{rk^2}_{\left(2K \sqrt{\gamma_1}  < (L_1+1) 
\epsilon \right) }(\Gamma_R(Y;m_1,k,\gamma_1)) 
\leq (L_1+1)^{rk^2} \cdot
H^{rk^2}_{\left( 2KL \sqrt{\gamma} < 
\epsilon \right)}(\Gamma_{LR}(Z;m,k,\gamma)).\]

\noindent By definition we now have 

\[ \overline{\mathbb H}^r_{2K, (L_1+1)\epsilon, R}(Y;m_1,\gamma_1) \leq
\overline{\mathbb H}^r_{2KL,\epsilon,LR}(Z;m,\gamma) + r \cdot \log (L_1
+1).\]

\noindent This hold for $m, m_1$ sufficiently large and $\gamma, \gamma_1$
sufficiently small.  Thus,

\[ \overline{\mathbb H}^r_{2K,(L_1+1)\epsilon,R}(Y) \leq \overline{\mathbb
H}^r_{2KL, \epsilon, L R}(Z) + r \cdot \log(L_1 +1).\]

\noindent Thus, for any $\epsilon >0$ we have $\overline{\mathbb
H}^r_{(L_1+ 1)\epsilon, R}(Y)  \leq \overline{\mathbb H}^r_{\epsilon,
LR}(Z) + r \cdot \log(L_1 +1).$ Taking a $\lim_{\epsilon \rightarrow 0}$
on both sides yields $\overline{\mathbb H}^r_R(Y) \leq \overline{\mathbb
H}^r_{LR}(Z) + r \cdot \log(L_1+1)$ where both $L$ and $L_1$ are dependent
upon $R.$

	Now suppose $\overline {\mathbb H}^r(Z) = -\infty.$ Then for all
$R>0$ $\overline{\mathbb H}^r_R(Z) = - \infty$ and by the last sentence of
the preceding paragraph this means $\overline{\mathbb H}^r_R(Y) = -
\infty$ for all such $R,$ whence $\overline{\mathbb H}^r(Y) = - \infty.$
By definition, $\overline{\mathbb H}(Y) \leq \overline{\mathbb H}(Z)$ as
desired.
\end{proof}

\section{The Free Hausdorff Dimension of Finite Dimensional Algebras}

In this section we show that if $M$ is finite dimensional and $\{z_1, 
\ldots,z_n\}$ generates $M$, then 
\[\mathbb H(z_1,\ldots,z_n) = \delta_0(z_1,\ldots,z_n).\] 

\noindent The argument is geometrically simple and it amounts to a
slightly finer analysis than that in [4] where the main objective was to
calculate $\delta_0(\cdot)$ of sets of self-adjoint elements which
generate hyperfinite von Neumann algebras.  The metric information of the
microstate space of $\{z_1,\ldots,z_n\}$ is encapsulated in the unitary
orbit of the images of $z_i$ under a single representation of $M$ which
preserves traces.  This set in turn, is Lipschitz isomorphic to the
homogeneous space obtained by quotienting $U_k$ by the subgroup of
consisting of those unitaries commuting with the image of $M$ under the
representation.  By [8] a neighborhood of this homogeneous space is
(modulo a Lipschitz distortion) a ball of operator norm radius $r$ in
Euclidean space of dimension $\alpha k^2$ where $\alpha =
\delta_0(z_1,\ldots,z_n)$ ($\alpha$ depends only on $M$ and $\varphi$).  
By the computations of [7] the asymptotic metric information of this set
is roughly that of a ball of $|\cdot|_2$-radius $r$ in Euclidean space of
dimension $\alpha k^2.$ The Hausdorff quantities of balls are easy to deal
with and yield the expected dimension.

Because such balls are $\alpha k^2$-sets in their ambient space of equal 
$|\cdot |_2$-radius, the 
argument above says that $\{z_1,\ldots,z_n\}$ is an $\alpha$-set 
and thereby shows that $\mathbb H(z_1,\ldots,z_n) = \delta_0(z_1,\ldots,z_n) = 
\alpha.$

We start with an upper bound which works for all hyperfinite von Neumann 
algebras and then proceed with the lower bound for the finite dimensional case.

\subsection{Upper Bound}

        Throughout the subsection assume $z_1,\ldots, z_n$ are selfadjoint
generators for $M$ and that $M$ is hyperfinite.  By decomposing $M$ over
its center

\begin{eqnarray}& & M \simeq M_0 \oplus (\oplus_{i=1}^{s}
M_{k_i}(\mathbb C)) \oplus M_{\infty} \nonumber \\ & &  \varphi \simeq
\alpha_0 \varphi_0 \oplus (\oplus_{i=1}^{s}
\alpha_i tr_{k_i}) \oplus 0 \nonumber
\end{eqnarray}

\noindent where $s \in \mathbb N \bigcup \{0\} \bigcup \{\infty\},
\alpha_i >0$ for $1 \leq i \leq s$ $( i \in \mathbb N),$ $M_0$ is a
diffuse von Neumann algebra or $\{0\}$, $\varphi_0$ is a faithful, tracial
state on $M_0$ and $\alpha_0 >0$ if $M_0 \neq \{0\},$ $\varphi_0 = 0$ and
$\alpha_0 =0$ if $M_0 = \{0\},$ and $M_{\infty}$ is a von Neumann algebra
or ${0}.$ Set $\alpha = 1 - \sum_{i=1}^s \frac {\alpha_i^2}{k_i^2}.$ We
show in this section that $\mathbb H^{\alpha}(z_1,\ldots,z_n)  < \infty.$

        We remark that by [4] $\delta_0(z_1,\ldots,z_n) = \alpha$ so by 
Corollary 3.5 $\mathbb H(z_1,\ldots,z_n) \leq \alpha.$  However, we have 
the slightly stronger statement below:

\begin{lemma} $\mathbb P^{\alpha}(z_1,\ldots,z_n) < \infty.$
\end{lemma}
\begin{proof} By Theorem 3.10 of [4] there exists a $C>0$ such 
that for $\epsilon >0$ sufficiently small

\[ \chi(a_1+ \epsilon s_1, \ldots, a_n + \epsilon s_n, I + \epsilon s_{n+1}: 
s_1,\ldots, s_{n+1}) \leq \log (\epsilon^{n+1 - \alpha}) + \log(4^{n+1}D_0) 
\]

\noindent where $D_0 = \pi^{n+1}(8(R+1))^{n+1}(C+1)6^{n+1}$ and $R$ is the
maximum of the operator norms of the $a_i.$ By [5] for $\epsilon >0$
sufficiently small

\begin{eqnarray*} \chi(a_1 + \epsilon s_1,\ldots, a_n + \epsilon s_n, I + 
\epsilon s_{n+1}; s_1,\ldots,s_{n+1}) & \geq &  \mathbb P_{2 \epsilon 
\sqrt{n}}(z_1,\ldots,z_n, I) \\ & + & (n+1) \cdot \log \epsilon + 
\chi(s_1,\ldots,s_{n+1}).
\end{eqnarray*}

\noindent  Hence

\[ \mathbb P_{2 \epsilon \sqrt{n}}(z_1, \ldots, z_n) \leq \mathbb P_{2 
\epsilon \sqrt{n}}(z_1, \ldots, z_n, I) \leq \log (\epsilon^{-\alpha}) + K \]

\noindent where $K = \log ((2\pi e)^{-\frac{n+1}{2}} 4^{n+1} D_0).$ By
Lemma 3.11 we're done. \end{proof}

More generally the analysis of [3] shows:

\begin{corollary} Suppose $\{z_1,\ldots,z_n\}$ generates $M$ and 
$\langle u_j \rangle_{j=1}^{s}$ 
is a sequence of Haar unitaries which also generates $M.$  If $u_{j+1}u_ju_{j+1}^* \in \{u_1,\ldots,u_j\}^{\prime \prime}$ for each $1 \leq j \leq s-1,$ then 

\[ \mathbb P^1(z_1,\ldots,z_n) < \infty. \]
\end{corollary}

\subsection{Lower Bound}

Throughout assume that $M = \oplus_{i=1}^p M_{k_i}(\mathbb C) \neq
\mathbb C I,$ $\varphi = \oplus_{i=1}^p \alpha_j tr_{k_i}$ where $p \in
\mathbb N$ and $\alpha_i >0$ for each $i,$ and $R> \max_{1 \leq j \leq
n}\{\|z_j\|\}.$ Also assume that the $z_j$ generate $M.$ Set $\alpha = 1 -
\sum_{i=1} \frac {\alpha_i^2}{k_i^2}.$ By Corollary 5.8 of [4] for any set
of self-adjoint generators $a_1,\ldots, a_m$ for $M,$
$\delta_0(a_1,\ldots, a_m) = \alpha.$

	We recall some basic facts from section 5 of [4].  There exists a
$z \in M$ such that the $*$-algebra $z$ generates is all of $M.$ For a
representation $\pi : M \rightarrow M_k(\mathbb C)$ define $H_{\pi}$ to be
the unitary group of $(\pi(M))^{\prime}$ and $X_{\pi} = U_k/H_{\pi}.$
Endow $X_{\pi}$ with the quotient metric $d_2$ derived from the $|
\cdot|_2$-metric on $U_k.$ Define $U_{\pi}(z)$ to be the unitary orbit of
$\pi(z).$ Consider the map $f_{\pi} : U_{\pi}(z)  \rightarrow X_{\pi}$
given by $f_{\pi}(u \pi(z) u^*) = q(u)$ where $q:  U_k \rightarrow
X_{\pi}$ is the quotient map.  By Lemma 5.4 of [4]

\[\{f_{\pi} : \hspace{.05in} \text {for some} \hspace {.08in} k \in 
\mathbb N \hspace{.05in} \pi : M \rightarrow M_k(\mathbb C) \hspace{.08in} 
\text {is a 
representation} \} \] 

\noindent has a uniform Lipschitz constant $D >1.$ Finally, there exists a
polynomial $f$ in $n$ noncommuting variables satisfying $f(z_1,\ldots,z_n)
= z.$ We find a constant $L >1$ such that for any $k\in \mathbb N$ and
$\xi_1, \ldots, \xi_n, \eta_1,\ldots,\eta_n \in (M_k^{sa}(\mathbb C))_R$

\[ |f(\xi_1,\ldots,\xi_n) - f(\eta_1,\ldots,\eta_n)|_2 \leq L \cdot \max 
\{|\xi_i - \eta_i|_2: 1 \leq i \leq n \}. \]  

	We will also need the following lemma.  It is a sharpening of
Lemma 3.6 of [4].
	
\begin{lemma}If $1>\varepsilon>0,$ then there exists a $N \in \mathbb N$ 
such that for any $k > N$ there is a corresponding 
$*$-homomorphism $\sigma_k:M \rightarrow M_k(\mathbb C)$ satisfying:

\begin{itemize} \item $\|tr_k \circ \sigma_k - \varphi \|< \varepsilon.$ 
\item $\dim (U_k/H_k) \geq \alpha k^2$ where $H_k$ is the unitary group of 
$\sigma_k(M)^{\prime}$ and $H_k$ is tractable in the sense defined in [4].
\end{itemize} \end{lemma}

\begin{proof} First suppose that for some $i,j$ $\alpha_i n_j^2 \neq
\alpha_j n_i^2.$ Without loss of generality we may assume that $\alpha_1
n_2^2 > \alpha_2 n_1^2.$ Given $\varepsilon >0$ as above we choose
$\varepsilon_1 < \varepsilon$ so that if $\beta_1 = \alpha_1 -
\varepsilon_1, \beta_2 = \alpha_2 +\varepsilon_1,$ and $\beta_i =
\alpha_i$ for $3 \leq i \leq p$ then $\sum_{i=1} \frac{\beta_i^2}{n_i^2} <
\sum_{i=1}^p \frac{\alpha_i^2}{n_i^2} - \delta$ for some $1> \varepsilon
>\delta >0$ dependent upon $\varepsilon_1.$ Obviously $\beta_1 + \cdots +
\beta_p =1.$

We now proceed as in Lemma 3.6.  Choose $n_0 \in \mathbb N$ such that
$\frac{1}{n_0} < \frac{\delta}{p^2}$ and set $k_0 = (n_0+1)n_1 \cdots
n_p.$ Suppose $k > k_0.$ Find the unique $n \in \mathbb N$ satisfying

\[ nn_1 \cdots n_p \leq k < (n+1)n_1 \cdots n_p.\]

\noindent Set $d = nn_1 \cdots n_p$ and find $m_1 \cdots m_p \in \mathbb N
\bigcup \{0\}$ satisfying $\beta_i - \frac{\delta}{4p^2} < 
\frac{m_i}{n} < \beta_i + \frac{\delta}{4p^2}$ and $\sum_{i=1}^p 
\frac{m_i}{n}=1.$ Set $l_i=\frac {d m_i}{n n_i} \in \mathbb N \bigcup \{0\}$ and $l_{p+1} = k -
\sum_{i=1}^p l_i n_i.$ Define $\sigma_k:N \rightarrow M_k(\mathbb C)$ by

\[\sigma_k(x_1,\ldots,x_n)= \begin{bmatrix} I_{l_1} \otimes x_1 & 0 
& \cdots & 0 \\ 0 & \ddots &  & \vdots \\ \vdots &  & I_{l_p}
\otimes x_p & 0 \\ 0 & \cdots & 0 & 0_{l_{p+1}} \\
\end{bmatrix} \]

\noindent where $0_{l_{p+1}}$ is the $l_{p+1} \times l_{p+1}$ 0
matrix and $I_{l_i} \otimes x_i$ is the $l_i n_i \times     
l_in_i$ matrix obtained by taking each entry of $x_i, (x_i)_{st}$, and
stretching it out into $(x_i)_{st} \cdot I_{l_i}$ where $I_{l_i}$ is
the $l_i \times l_i$ identity matrix. 

\[ (tr_k \circ \sigma_k)(x_1,\ldots,x_p) = \frac {1}{k} \cdot \sum_{i=1}^p
l_i  \cdot Tr(x_i) = \sum_{i=1}^{p} \frac {d m_i}{k n } \cdot
tr_{n_i}(x_i). \] $\frac {d}{k} > 1- \frac{\delta}{p^2}$ so $\alpha_i +
\frac{\delta}{p^2} \geq \frac {d}{k} \cdot \frac {m_i}{n} > (\alpha_i -
\frac{\delta}{p^2})(1- \frac{\delta}{p^2}) > \alpha_i - 
\frac {\varepsilon}{p}.$  It follows that
$\|tr_k \circ \sigma_k - \varphi\| < \varepsilon $.

      $H_k$ consists of all matrices of the form \[ \begin{bmatrix}
u_1\otimes I_{n_1} & 0 & \cdots & 0 \\ 0 & \ddots & & \vdots \\ \vdots
& & u_p \otimes I_{n_p} & 0 \\ 0 & \cdots & 0 & u_{p+1} \\
\end{bmatrix} \] where $u_i \in U_{l_i}$ for $1 \leq i \leq p+1$ and $u_i
\otimes I_{n_i}$ is the $l_i n_i \times l_i n_i$ matrix obtained by
repeating $u_i \hspace{.1in} n_i$ times along the diagonal. $H_k$ is
obviously tractable.  Thus we have the estimate: 

\[ l_{p+1} = k -
\sum_{i=1}^{p} \frac{d m_i}{n} = k - d < n_1 \cdots n_p. \] 

\noindent So \begin{eqnarray*}\dim H_k = l_{p+1}^2 + \sum_{i=1}^{p}l_i^{2}
< k^2 \left (\frac{n_1^2\cdots n_p^2}{k^2} + \sum_{i=1}^{p} \frac
{m_i^2}{n^2 n_i^2} \right) & < & k^2 \left ( \frac{n_1^2\cdots n_p^2}{k^2}
+ \frac{\delta}{2p} + \sum_{i=1}^p \frac{\beta_i^2}{n_i^2} \right) \\ & <
& k^2 \left ( \frac{n_1^2\cdots n_p^2}{k^2} - \frac {\delta}{2} +
\sum_{i=1}^p \frac{\alpha_i^2}{n_i^2} \right) \\ & < & k^2 \cdot \left( 
\sum_{i=1}^p \frac{\alpha_i^2}{n_i^2} \right).
\\ \end{eqnarray*}

\noindent Hence, $\dim(U_k/H_k) \geq \alpha k^2.$ For $k > k_0 = N$ we
have produced a $*$-homomorphism $\sigma_k:M \rightarrow M_k(\mathbb C)$
satisfying all the properties of the lemma.

Now suppose that $\alpha_i n_j^2 = \alpha_j n_i^2$ for all $i,j.$ It
follows that $\alpha_i \in \mathbb Q$ for all $i.$ Otherwise,
$\alpha_i$ is irrational for some $i.$ Thus, $\sum_{i=1}^p n_j^2/n_i^2 =
\alpha_i^{-1}$ which is absurd.  For each $i$ write $\alpha_i = p_i/q_i,$
$p_i, q_i \in \mathbb N.$ Set $N = q_1 \cdots q_p n_1 \cdots n_p.$ Suppose
$k = k_1 N$ for some $k_1 \in \mathbb N.$ Define $l_i = \frac{k
\alpha_i}{n_i} \in \mathbb N.$ As in the preceding argument define
$\sigma_k:N \rightarrow M_k(\mathbb C)$ by

\[\sigma_k(x_1,\ldots,x_n)= \begin{bmatrix} I_{l_1} \otimes x_1 &  
  & 0\\  & \ddots & \\ 0 &  & I_{l_p}
\otimes x_p  \\
\end{bmatrix}. \]

\noindent It is plainly seen that $tr_k \circ \sigma_k = \varphi$ and that
$H_k$ is a tractable subgroup with $\dim H_k = k^2 \cdot \sum_{i=1}^p
\frac{\alpha_i^2}{n_i^2}.$ Thus $\dim (U_k/H_k) = \alpha k^2$ and we have
the desired result for all multiples $k$ of $N.$ It is easy from here to
show that the result holds for all sufficiently large $k$ and we leave the
proof to the reader. \end{proof} \vspace{.1in}

\begin{lemma} $\{z_1,\ldots,z_n\}$ is an $\alpha$-set. \end{lemma}

\begin{proof} By Lemma 4.1 it suffices to show that $\mathbb
H^{\alpha}(z_1,\ldots,z_n) > -\infty.$ Recall the proof of Lemma 5.2 in [4].  
Replacing Lemma 3.6 of [4] in the proof with Lemma 4.1 above, the
arguments of [4] produce $1 > \lambda, \zeta, r >0$ such that for any
given $m \in \mathbb N$ and $\gamma >0$ there exists an $N \in \mathbb N$ 
such that for $k \geq N$ there exists a 
$*$-homomorphism $\sigma_k : M \rightarrow M_{k}(\mathbb C)$ and:

\begin{itemize}

\item $\|tr_k \circ \sigma_k - \varphi \| 
< \frac {\gamma}{(R+1)^m}.$

\item The set of unitaries $H_k$
of $\sigma_k(M)^{\prime}$ is a tractable Lie subgroup of $U_k$ and setting
$X_k = U_k/H_k,$ $ \dim(X_k) \geq \alpha k^2.$

\item  Define $\mathcal H_k 
\subset iM^{sa}_k(\mathbb C)$ to be the Lie 
subalgebra of $H_k$ (as above) and $\mathcal X_k$ to be the orthogonal 
complement of 
$\mathcal H_k$ with respect to the Hilbert-Schmidt inner product.  For 
every $s>0$ write $\mathcal X^s_k$ for the ball in $\mathcal X_k$ of 
operator norm less than or equal to $s$ and $c_k$ for the volume of the 
ball of $\mathcal X_k$ of $|\cdot|_2$ of radius 1.  Here all volume 
quantities are obtained from Lebesgue measure when the spaces are 
given the real inner product induced by $Tr.$  

\[ \frac{\text{vol}(\mathcal X^1_k)}{c_k} > (\zeta)^{\dim \mathcal X_k}.\]

\item For any $x,y \in \mathcal X^r_k$
\[d_2(q(e^x), q(e^y)) \geq \lambda |x-y|_2.\]

\end{itemize}

\noindent Suppose $m$ and $\gamma$ are fixed and $k$ so that the four
conditions above hold.  Suppose also that $\epsilon < \lambda ( DL)^{-1}.$
If $T_k$ denotes the unitary orbit $\{(u\sigma_k(z_1)u^*,\ldots,
u\sigma_k(z_n)u^*) : u \in U_k\},$ then clearly

\[\mathbb H^{\alpha}_{\epsilon, R}(z_1,\ldots,z_n;m,\gamma) \geq 
\limsup_{k \rightarrow \infty} \left (k^{-2} \cdot \log H^{\alpha 
k^2}_\epsilon(T_k)\right).\]

Define $g_k: X_k \rightarrow T_k$ by $g_k(q(u))= (u\sigma_k(z_1) u^*,
\ldots, u\sigma_k(z_n)u^*).$ Denote by $\Omega_k$ the image of $\mathcal
X^r_k$ under the map $g_k \circ q \circ e.$ Clearly $\Omega_k \subset T_k$
and the map $\Phi = (q \circ e)^{-1} \circ f_{\pi} \circ f : \Omega_k
\rightarrow \mathcal X^r_k$ is a well-defined (by the fourth condition
above) surjective map with $\| \Phi \|_{Lip} \leq \frac {DL}{\lambda}.$
Hence

\[H^{\alpha k^2}_{\epsilon}(T_k) \geq H^{\alpha k^2}_{\epsilon}(\Omega_k)  
\geq \left(\frac{\lambda}{DL}\right )^{\alpha k^2} \cdot H^{\alpha
k^2}_{DL \epsilon \lambda^{-1}}(\mathcal X^r_k).\]

\noindent Suppose $\langle \theta_j \rangle_{j \in J}$ is  
a countable cover of $\mathcal X^r_k.$  We have by volume comparison

\[ \sum_{j\in J} \frac{(\sqrt{\pi k}|\theta_j|)^{\dim \mathcal
X_k}}{\Gamma(\frac{\dim \mathcal X_k}{2} +1)} \geq \sum_{j \in J}
\text{vol}(\theta_j) \geq \text{vol}(\mathcal X^r_k) \geq (r \zeta)^{\dim 
\mathcal X_k} \cdot \frac {(\sqrt{\pi k})^{\dim \mathcal 
X_k}}{\Gamma(\frac{\dim \mathcal X_k}{2} +1)}.\]

\noindent Thus 

\[ H^{\alpha k^2}_{DL \epsilon \lambda^{-1}}(\mathcal
X_k^r) \geq H^{\dim X_k}_{DL \epsilon \lambda^{-1}}(\mathcal X_k^r) \geq
(r \zeta)^{\dim \mathcal X_k}. \]

\noindent Following the chains of inequalities for such $k \geq N$ \[
H_{\epsilon}^{\alpha k^2}(T_k) \geq \left( \frac{\lambda r
\zeta}{DL}\right)^{\dim X_k}. \]

\noindent It follows that

\[ \limsup_{k \rightarrow \infty} \left (k^{-2}\cdot \log 
H^{\alpha k^2}_{\epsilon} 
(T_k) \right) \geq  \log 
\left (\frac{\lambda r \zeta }{DL} \right) .\]

\noindent By the concluding inequality of the preceding paragraph $\mathbb
H^{\alpha}_{\epsilon, R}(z_1,\ldots,z_n;m, \gamma)$ exceeds the right hand
expression above.  Forcing $\epsilon \rightarrow 0$ we conclude that
$\mathbb H^{\alpha}_R (z_1,\ldots,z_n)$ exceeds the right hand expression
above. $\mathbb H^{\alpha} (z_1,\ldots,z_n) > -\infty.$
\end{proof} 

\begin{corollary} $\mathbb H(z_1,\ldots,z_n) = \alpha =
\delta_0(z_1,\ldots,z_n).$ \end{corollary}

\section{The free Hausdorff dimension of a single selfadjoint }

	We show that the free Hausdorff dimension and modified free
entropy dimension are equal for single selfadjoints.  In the first
subsection we prove an easy lemma on finding lower bounds for Hausdorff
measure quantities of locally isometric spaces.  From there we compute in
the second part their asymptotic limit to arrive at the desired claim.  
Finally, we present a microstates version of the classical fact that a
space with Hausdorff dimension strictly less than 1 is totally
disconnected.

\subsection {A Lemma on Hausdorff Measures}

	Finding sharp lower bounds for the free Hausdorff dimension of a
given $n$-tuple hinges on estimating $H^{\alpha k^{2}}_{\epsilon}$ of the
microstate spaces.  Here $\alpha$ and $\epsilon$ remain fixed as $k$ tends
to infinity.  The lemma we will prove below says that for locally
isometric spaces (metric spaces such that any two $\epsilon$ balls are
isometric), the right lower bounds on the $\epsilon$ packing numbers give
the right lower bounds on $H^{\alpha k^2}_{\epsilon}.$ We use this result
in the next subsection through the following argument.  The microstate
spaces of a single self-adjoint are unitary orbits of single self-adjoint
matrices with appropriate eigenvalue densites.  Such sets are locally
isometric and the volumes of the $\epsilon$-neighborhoods of such orbits
are well known [6].  Invoking the lemma below with appropriate bounds will
then provide the result.

\begin{lemma}Suppose $X \subset \mathbb R^d$ is a $g$-set such that any
two open $\epsilon$ balls of $X$ are isometric.  Assume further that for
some $\alpha \geq 0$ there exist $C > 1 >\epsilon_0>0,$ such that for any
$\epsilon_0 > \epsilon >0$

\[ P_{\epsilon}(X) > C \cdot \epsilon^{-\alpha}.\]

\noindent Then for any Borel set $E \subset X$

\[ H^\alpha_{\epsilon_0}(E)  >  C \cdot 
\frac{H^g(E)}{H^g(X)}.
\]
\end{lemma}

\begin{proof} Pick a point $x \in X$ and for $\epsilon >0$ denote by
$D_{\epsilon}$ the open ball in $X$ of radius $\epsilon$ centered at
$x.$ Using the lower bound on $P_{\epsilon}(X)$ and the fact that any two
open $\epsilon$ balls of $X$ are isometric we have for any $\epsilon_0 > 
\epsilon >0$
 
\[ H^{g}(D_{\epsilon}) \leq \frac {H^g(X)}{P_{\epsilon}(X)} \leq 
\frac{H^g(X)}{C} \cdot \epsilon^{\alpha}.\]

\noindent Now suppose $\langle \theta_j \rangle_{j \in J}$ is an
$\epsilon_0$-cover for $E.$ For $\delta > 0$ and for each $j$ find an open
ball $D_j$ of $X$ with radius no greater than $(1 + \delta)|\theta_j|$
satisfying $\theta_j \subset D_j.$ Using the inequality above paired with
the assumption that any two open balls of $X$ of equal radius are
isometric we now have the estimate:

\[ (1+ \delta)^{\alpha} \cdot \sum_{j \in J} |\theta_j|^{\alpha} \geq
\frac {C}{H^g(X)} \cdot \sum_{j \in J} H^g(D_j) \geq C \cdot
\frac{H^g(E)}{H^g(X)}. \]

\noindent Forcing $\delta \rightarrow 0$ it follows that $\sum_{j \in J}
|\theta_j|^{\alpha} \geq C \cdot \frac{H^g(E)}{H^g(X)}.$ $\langle \theta_j
\rangle_{j \in J}$ being an arbitrary $\epsilon_0$-cover for $E$ the
conclusion follows. \end{proof}

\subsection{The Estimates}

It's now just a matter of putting the lemma together with some strong
packing estimates.

\begin{lemma} If $z = z^* \in M,$ then $\mathbb H(z) = \delta_0 (z).$	  
\end{lemma}

\begin{proof} By Corollary 3.5 $\mathbb H(z) \leq \delta_0(z)$ so it
suffices to prove the reverse inequality.  Recall from [10] that
$\delta_0(z) = 1 - \sum_{t \in sp(z)} \mu(\{t\})^2,$ where $\mu$ is the
Borel measure induced on $sp(z)$ by $\varphi,$ so we have to show that
$\mathbb H(z) \geq 1 - \sum_{t \in sp(z)} \mu(\{t\})^2.$ 

Write $\mu = \sigma + \nu$ where $\sigma$ is the atomic part of $\mu$ and
$\nu$ is the diffuse part of $\mu.$ $\sigma = \sum_{i=1}^s c_i
\delta_{r_i}$ for some $s \in \mathbb N \bigcup \{0\} \bigcup \{\infty\},$
$c_i \geq c_{i+1} > 0,$ and where for $i \neq j,$ $r_i \neq r_j.$ Suppose
$R > \|z\|$ and $\frac{1}{2} >\tau >0.$ We can find $l$ such that $
\sum_{i=l+1}^s c_i < \frac{\tau}{3}.$ For $\epsilon >0$ define
$D_{\epsilon} = \{(s,t) \in [a,b]^2: |s-t| < \epsilon \}.$ Because $\nu$
is diffuse by Fubini's theorem for small enough $\epsilon_0 <1,$

\[(\nu \times\nu)(D_{3 \epsilon_0}) < \tau.\]

\noindent Arrange it so that $\epsilon_0$ also satisfies $\nu([r_i -
\epsilon_0, r_i+ \epsilon_0]) < \frac{\tau}{3l}$ and $|r_i - r_j| > 3
\epsilon_0$ for $1 \leq i < j \leq l.$ Set $c = \nu([a,b]).$ Again because 
$\nu$ is diffuse for each $k$ and $1 \leq j \leq
[ck]$ there exists a largest number $\lambda_{jk} \in [a, b]$ satisfying
$\nu([a, \lambda_{jk}])= \frac{j}{k}.$

Now observe that for $k$ large enough there exists a subset $S_k$ of $\{1,
\ldots, [ck]\}$ such that for any $i \in S_k, |\lambda_{ik} -
\lambda_{(i+1)k}| < \epsilon_0$ and $\#S_k > (1 - \frac{\tau}{3})[ck].$ To
see this for each $k$ consider the maximum number of elements in a subset
of $\{1,\ldots, [ck]\}$ such that for any element $i$ in this subset,
$|\lambda_{ik} - \lambda_{(i+1)k}| < \epsilon_0.$ Find a subset $S_k$
which achieves this maximum number.  For any $i \notin S_k,$
$|\lambda_{ik} -\lambda_{(i+1)k}| \geq \epsilon_0.$ Define $E_k = \{1,
\ldots, [ck]\} - S_k$ and denoting by $m$ Lebesgue measure on $\mathbb R$

\[ b-a > m(\cup_{i \in E_k} [\lambda_{ik}, \lambda_{(i+1)k}]) \geq \# E_k
\cdot \epsilon_0. \]

\noindent Consequently $\# E_k < \frac{b-a}{\epsilon_0}$ and thus for $k$
large enough $\# E_k < \frac{\tau [ck]}{3}$ (provided $c >0;$ if $c=0,$ 
then the claim is vacuously satisfied), whence $\#S_k \geq (1-
\frac{\tau}{3})[ck].$

By what has preceded for $k$ sufficiently large I can find a subset $G_k$
of $\{\lambda_{1k},\ldots,\lambda_{[ck]k} \}$ such that

\begin{itemize}

\item For any $i \in G_k,$ $|\lambda_{ik} - \lambda_{(i+1)k}| < 
\epsilon_0.$ 

\item Any element
of $G_k$ is at least $3 \epsilon_0$ apart from any element of $\{r_1,
\ldots, r_l\}.$ 

\item $\#G_k > [ck] - \frac{\tau [ck]}{3} - \frac{\tau k}{3}.$ 

\end{itemize}

For each $k$ define $p_k = \#G_k + \sum_{i=1}^l [c_ik] ( > k - \tau k)$
and $A_k$ to be the diagonal $ p_k \times p_k$ matrix obtained by filling
in the the first $\#G_k$ diagonal entries with the elements of $G_k$
(ordered from least to greatest) and the last $\sum_{j=1}^l [c_jk]$
diagonal entries filled with $r_1$ repeated $[c_1k]$ times, $r_2$ repeated
$[c_2k]$ times, etc., in that order.  Define $B_k$ to be the $ (k - p_k )
\times (k - p_k)$ diagonal matrix obtained by filling in the entries
(ordered from least to greatest) with $\{\lambda_{1k}, \ldots,
\lambda_{[ck]k}\}- G_k,$ the entries $r_i$ repeated $[c_i k]$ for $i >l,$
and $0$ repeated as many times as necessary to make $B_k$ a $(k-p_k)
\times (k- p_k)$ matrix. Finally define $y_k$ to be the $k \times k$
matrix

\[
\begin{bmatrix}
A_k & 0 \\
0 &  B_k \\
\end{bmatrix}
\]

\noindent For any $m \in \mathbb N$ and $\gamma 
>0$ given, $y_k \in \Gamma_R(z;m,k,\gamma)$ for sufficiently large $k.$
For any $x \in M^{sa}_k(\mathbb C)$ denote by $\Theta(x)$ the unitary 
orbit of $x.$  We have that $\Theta(y_k) \in \Gamma_R(z;m,k,\gamma)$ for 
sufficiently large $k.$  

Now we want strong lower bounds for the packing numbers of 
$\Gamma_R(z;m,k,\gamma).$  Such bounds for $\Theta(y_k)$ will suffice. 
Write $\Theta_r(y_k)$ for the set of all $k \times k$ matrices of the form 

\[
\begin{bmatrix}
uA_ku^* & 0 \\
0 &  B_k \\
\end{bmatrix}
\]

\noindent where $u$ is a $p_k \times p_k$ unitary.  $\Theta_r(y_k) \subset 
\Theta(y_k)$ and because $1/2 > \tau$ 

\[ P_{\epsilon}(\Theta(y_k)) \geq 
P_{\epsilon}(\Theta_r(y_k)) \geq  P_{2 \epsilon}(\Theta(A_k)).\]

\noindent If we can find strong lower bounds for the packing numbers of
$\Theta(A_k),$ or equivalently strong lower bounds for the volume of the
$\epsilon$-neighborhoods of $\Theta(A_k),$ then we can invoke Lemma 6.5
and arrive at a lower bound for the Hausdorff quantities of
$\Gamma_R(z;m,k,\gamma).$

Denote $G$ to be the group of diagonal unitaries and $\mathbb R^k_{<}$ to
be the set of all $(t_1,\ldots,t_k) \in \mathbb R^k$ such that $t_1 <
\cdots < t_k.$ There exists a map $ \Phi : M^{sa}_k(\mathbb C) \rightarrow
U_k/G \times \mathbb R^k_{<}$ defined almost everywhere on
$M^{sa}_k(\mathbb C)$ such that for each $x \in M^{sa}_k(\mathbb C)$
$\Phi(x) = (h, z)$ where $z$ is a diagonal matrix with real entries
satisfying $z_{11} < \cdots < z_{kk}$ and $h$ is the image of any unitary
$u$ in $U_k/G$ satisfying $uzu^* = x.$ By results of Mehta [6] the map
$\Phi$ induces a measure $\mu$ on $U_k/G \times \mathbb R^k_<$ given
by $\mu(E) = \text{vol} (\Phi^{-1}(E))$ and moreover,

\[ \mu = \nu \times D_k \cdot
\int_{\mathbb R^k_{<}}
\Pi_{i < j} (t_i - t_j)^2 \, dt_1 \cdots dt_k \]

\noindent where $D_k =\frac{\pi^{k(k-1)/2}}{\Pi_{j=1}^k j!}$ and $\nu$ is
the normalized measure on $U_k/G$ induced by Haar measure on $U_k.$ Write
$\Theta_{\epsilon}(A_k)$ for the $|\cdot|_2$
$\epsilon$-neighborhood of the unitary orbit of $A_k$ and $\Theta(A_k)$ 
for the unitary orbit of $A_k.$  A matrix will be 
in $\Theta_{\epsilon}(A_k)$ iff the sequence obtained by listing its
eigenvalues in increasing order and according to multiplicity, differs
from the similar sequence obtained from the eigenvalues of $y_k$ by no 
more than $\sqrt{p_k} \cdot \epsilon$ in $\ell^2$ norm.  In particular 
this will happen if the $j$th terms of the sequences differ by no more
than $\epsilon.$  

Now for each $k$ write $a_{1k}, \ldots, a_{p_k k}$ for the eigenvalues
of $A_k$ ordered from least to greatest and according to multiplicity.
Consider the region in $\mathbb R^{p_k}$ obtained by taking the Cartesian
product

\[ [a_{1k} - \epsilon, a_{1k} + \epsilon] \times \cdots \times    
[a_{p_k k} - \epsilon, a_{p_k k} + \epsilon]
\]

\noindent Denote by $\Omega_k$ the intersection of this region with
$\mathbb R^{p_k}_{<}.$ Integrating over $\Omega_k$ according to the
density given above it follows that for $\epsilon_0 > \epsilon >0$
$\text{vol}(\Theta_{\epsilon}(A_k))$ exceeds

\begin{eqnarray} D_{p_k} \cdot \int_{\Omega_k} \Pi_{i <j}(t_i -t_j)^2 \,
dt_1
\cdots dt_{p_k}. \end{eqnarray}

\noindent Denote by $W_k$ all 2-tuples $(i,j)$ such that $1 \leq i < j
\leq p_k$ and $|a_{ik} -a_{jk}| < \epsilon_0.$ Generously estimating, 
(1) dominates

\begin{eqnarray} D_{p_k} \cdot \epsilon_0^{k^2} \cdot \int_{\Omega_k}
\Pi_{(i,j) \in
W_k} (t_i - t_j)^2 \, dt_1 \cdots dt_{p_k}.  
\end{eqnarray}

\noindent Consider the map $F: [- \epsilon, \epsilon]^{p_k} \cap \mathbb
R^{p_k}_{<} \rightarrow \Omega_k$ which sends $(t_1,\ldots, t_{p_k})$ to
$(a_{1k} + t_1, \ldots, a_{p_k k} + t_{p_k}).$  By a change of variables
formula via this map (2) dominates

\begin{eqnarray} D_{p_k} \cdot \epsilon_0^{k^2}\cdot \int_{[- \epsilon,
\epsilon]^{p_k} \cap \mathbb R^{p_k}_<} \Pi_{(i,j) \in W_k} (t_i - t_j)^2
\, dt_1 \cdots dt_{p_k}. \end{eqnarray}

\noindent By Selberg's integral formula and a change of variables we have

\begin{eqnarray*} \epsilon^{p_k^2} \cdot \Pi_{j=1}^{p_k} \frac{\Gamma(j+2)
\Gamma(j+1)^2}{\Gamma(p_k+j+1)} & = & \int_{[-\epsilon, \epsilon]^{p_k}}
\Pi_{ i < j} (t_i -t_j)^2 \, dt_1 \cdots dt_{p_k} \\ & = & p_k! \cdot
\int_{[-\epsilon, \epsilon]^{p_k} \cap \mathbb R^{p_k}_{<}} \Pi_{i < j}
(t_i - t_j)^2 \, dt_1 \cdots dt_{p_k} \\ & < & p_k! \cdot
(2\epsilon)^{p_k^2 - \#W_k} \cdot \int_{[- \epsilon, \epsilon]^{p_k} \cap
\mathbb R^{p_k}_< } \Pi_{(i,j) \in W_k} (t_i - t_j)^2 \, dt_1 \cdots
dt_{p_k} . \end{eqnarray*}

\noindent Thus, 

\[ (p_k!)^{-1} \cdot 2^{-k^2} \epsilon^{\# W_k} \cdot
\Pi_{j=1}^{p_k} \frac{\Gamma(j+2) \Gamma(j+1)^2}{\Gamma(p_k+j+1)} 
< \int_{[- \epsilon, \epsilon]^{p_k} \cap \mathbb R^{p_k}_< }
\Pi_{(i,j) \in W_k} (t_i - t_j)^2 \, dt_1 \cdots dt_{p_k}.  \]

All I need to do now is find an upper bound $\#W_k.$ Write $W_k = T_k
\cup V_k$ where $T_k$ consists of all $(i,j)$ such that $a_{ik} = a_{jk} =
r_l$ for some $1 \leq l \leq n$ and $V_k = W_k - T_k.$ Because $|r_i -
r_j| > 3 \epsilon_0$ for $1 \leq i < j \leq l,$ it follows that $\# T_k =
\sum_{j=1}^l [c_j k]^2.$ Let's estimate $V_k.$ First observe that if
$(i,j) \in V_k,$ then either $a_{ik}$ or $a_{jk}$ does not lie in
$\{r_1,\ldots, r_l\};$ consequently, both $a_{ik}$ and $a_{jk}$ are not in
$\{r_1,\ldots, r_l\}$ because all elements of $G_k$ are at least $3
\epsilon_0$ away from $r_1,\ldots, r_l.$ Therefore they are elements of
$G_k.$ For each $(i,j) \in V_k$ denote by $S(i,j)$ the closed square
$[a_{ik}, a_{(i+1)k}] \times [a_{jk}, a_{(j-1)k}].$ Because $|a_{ik} -
a_{(i+1)k}|, |a_{jk} - a_{(j+1)k}| < \epsilon_0,$ for each $(i,j) \in V_k,
S(i,j) \in D_{3 \epsilon_0}.$ Also, for $(i^{\prime}, j^{\prime}) \in V_k,
S(i,j) \cap S(i^{\prime}, j^{\prime})$ has $\nu \times \nu$ measure $0$
because $\nu$ is diffuse.

\[ \tau > (\nu \times \nu) (D_{3 \epsilon_0}) \geq (\nu \times \nu)
(\cup_{(i,j) \in W_k} S(i,j)) = \sum_{(i,j) \in W_k} (\nu \times \nu)
(S(i,j)) = \#V_k \cdot k^2.\]

\noindent Consequently, $\#W_k < \#T_k + \#V_k < (\tau + \sum_{j=1}^l
c_j^2)k^2.$ Write $\beta = \tau + \sum_{j=1}^l c_j^2.$ Substituting this
into the previous inequality we now have

\[ (p_k!)^{-1} \cdot 2^{-k^2} \cdot \epsilon^{\beta k^2} \cdot 
\Pi_{j=1}^{p_k}
\frac{\Gamma(j+2) \Gamma(j+1)^2}{\Gamma(k+j+1)} < \int_{[- \epsilon,
\epsilon]^{p_k} \cap \mathbb R^{p_k}_< } \Pi_{(i,j) \in W_k} (t_i - t_j)^2
\, dt_1 \cdots dt_{p_k}.\]

\noindent It follows that (3) dominates 

\[ D_{p_k} \cdot \epsilon_0^{k^2} \cdot (p_k!)^{-1} \cdot 2^{-k^2} \cdot
\epsilon^{\beta k^2} \cdot \Pi_{j=1}^{p_k} \frac{\Gamma(j+2)
\Gamma(j+1)^2}{\Gamma(k+j+1)} \]

\noindent and because $(1) > (2) > (3)$ in the previous paragraph, we
have that for $\epsilon_0 > \epsilon >0,
\text{vol}(\Theta_{\epsilon}(A_k)) >
L_k \cdot \epsilon^{\beta k^2}$ where 

\[L_k = D_{p_k} \cdot \epsilon_0^{k^2} \cdot (p_k!)^{-1} \cdot 2^{-k^2}
\cdot \Pi_{j=1}^{p_k^2} \frac{\Gamma(j+2)\Gamma(j+1)^2}{\Gamma(k+j+1)}.  
\]

For $\epsilon_0 > \epsilon>0$

\[P_{\epsilon}(\Theta(y_k)) \geq P_{2 \epsilon}(\Theta(A_k))  \geq L_k
\cdot \frac{\epsilon^{\beta k^2  - p_k^2} \Gamma(\frac{p_k^2}{2} +1 ) 
}{4^{k^2} (
\pi k)^{\frac{p_k^2}{2}} }.\]

\noindent Notice that $L_k$ is independent of
$\epsilon.$ By Lemma 6.1 it follows that for each $m \in \mathbb N$ and
$\gamma>0,$ and for sufficiently large $k$

\[ H^{p_k^2 - \beta k^2}_{\epsilon_0}(\Gamma_R(z;m,k,\gamma)) \geq 
\frac{L_k \Gamma(\frac{p^2_k}{2}+1)} {4^{k^2} \pi^{\frac 
{p_k^2}{2}} \sqrt{k}^{p_k^2} }.\]

\noindent Now $p^2_k - \beta k^2 > \alpha k^2$ where $\alpha = 1 -
\sum_{j=1}^l c_j^2 -4 \tau$ so for any $m \in \mathbb N$ and $\gamma >0$
$\mathbb H^{\alpha}_{\epsilon_0, R}(z;m,\gamma)$ dominates (by Stirling's
formula)

\begin{eqnarray*} && \liminf_{k\rightarrow \infty} \left( k^{-2} \log (L_k
\cdot \Gamma(\frac{p_k^2}{2}+1)) - \frac{p_k^2}{2k^2} \cdot \log k - 16\pi
\right) \\ & > & \liminf_{k \rightarrow \infty} \left (-k^{-2} \cdot \log
\Pi_{j=1}^{p_k} j! + \frac{p_k^2}{2} \log p_k \right ) + \log \epsilon_0
\\ & & + \liminf_{k \rightarrow \infty} k^{-2} \cdot \log \left
(\Pi_{j=1}^{p_k} \frac{\Gamma(j +2) \Gamma(j+1)^2}{\Gamma(k + j +1)}
\right) -16 \pi \\ & > & \log \epsilon_0 - 17 \pi. \end{eqnarray*}

\noindent $\liminf_{k \rightarrow \infty} k^{-2} \log \left
(\Pi_{j=1}^{p_k} \frac{\Gamma(j+2) \Gamma(j+1)^2}{\Gamma(k + j +1)}
\right)  > - \frac{5}{4}$ and $\liminf_{k \rightarrow \infty} \left( -
k^{-2} \log \Pi_{j=1}^{p_k} j! + \frac {p^2_k}{2 k^2} \log p_k \right)
>\frac {1}{4}$ above.  Both of these inequalities can be obtained from
some calculus. This lower bound being uniform in $m$ and $\gamma$

\[ \mathbb H^{\alpha}(z) \geq \mathbb
H^{\alpha}_{\epsilon_0,R}(z;m,\gamma)  > \log \epsilon_0 - 17 \pi.\]

\noindent $\alpha = 1 - 4\tau - \sum_{j=1}^l c_j^2.$ Finally since
$\frac{1}{2} > \tau >0 $ was arbitrary and $l \rightarrow s$ and $\alpha
\rightarrow 1 - \sum_{j=1}^s c_j^2 = \delta_0(z)$ (from [9] and [10]) as $
\tau \rightarrow 0$ it follows that $\mathbb H^{\alpha}(z) > - \infty$
for all $\alpha < \delta_0(z).$ $\mathbb H(z) \geq \delta_0(z).$
\end{proof}

\subsection{Minimal Projections}

From [4] if $\delta_0(z_1, \ldots, z_n) <1,$ then $\{z_1,\ldots,
z_n\}^{\prime \prime}$ has a minimal projection.  We will end this section
by showing that the same holds if $\delta_0$ is replaced by $\mathbb H.$
This is a slightly stronger statement since $\delta_0$ dominates $\mathbb
H.$ The corresponding classical fact is that a metric space with Hausdorff
dimension strictly less than 1 must be totally disconnected.  We will more
or less proceed by using the same argument in [4] for "weak hyperfinite
monotonicity of $\delta_0$" except we will limit our case to the situation
where the hyperfinite subalgebra is commutative.  Even then, the argument,
though intuively simple, requires more care since the quantities involved
are Hausdorff ones, and thus, harder to bound from below than the packing
quantities of $\delta_0.$ First a simple corollary from the computation
we've made.

\begin{corollary} If $z_k \rightarrow z$ strongly, then 
$\liminf_{k \rightarrow \infty} \mathbb H(z_k) \geq \mathbb H(z).$
\end{corollary}

\begin{proof} By Lemma 6.2 $\mathbb H(z) = \delta_0(z)$ and for 
all $k, \mathbb H(z_k) = \delta_0(z_k).$  Thus, by [10], $\liminf_{k 
\rightarrow \infty} \mathbb H(z_k) = \liminf_{k    
\rightarrow \infty} \delta_0(z_k) \geq \delta_0(z) = \mathbb H(z).$
\end{proof}

\begin{lemma} If $\{z_1, \ldots, z_n\}$ has finite dimensional 
approximants and $z=z^* \in \{z_1,\ldots, z_n\}^{\prime 
\prime},$ then  
\[\mathbb H(z_1,\ldots, z_n) \geq \mathbb H(z).\]
\end{lemma}

\begin{proof} We will first prove the statement under the additional
assumption that $z$ lies in the algebra $A$ generated by $\{z_1,\ldots,
z_n\}.$ Under this assumption there exists a polynomial $f$ in $n$
noncommuting variables, such that $f(z_1,\ldots, z_n) =z$ and we can
also assume that for any $n$ selfadjoint operators $h_1,\ldots, h_n$ on
a Hilbert space $f(h_1,\ldots, h_n)$ is again selfadjoint.  Now suppose
$R >0$ exceeds the operator norms of the $z_i.$ There exists an $L >
\|z\|$ such that for any $(h_1, \ldots, h_n ) \in (M^{sa}_k(\mathbb
C)_R)^n, f(h_1, \ldots, h_n) \in M^{sa}_k(\mathbb C)_L.$ Suppose
$\frac{1}{2} > \tau >0,$ and consider all the associated quantities
defined with respect to this $\tau$ in Lemma 6.2 for $z =
z^*: l, \epsilon_0, A_k, B_k, W_k, y_k, p_k, D_{p_k}.$

Suppose $m \in \mathbb N$ and $\gamma >0.$ By Lemma 4.2 of [4] there exist
$m_1 \in \mathbb N,$ and $\gamma_1 >0$ such that for any $a, b \in
\Gamma_L(z;m,k,\gamma),$ there exists a $u \in U_k$ satisfying $|uau^*
-b|_2 < \epsilon_0 \tau.$ We can choose $m_2 \in \mathbb N, m_2 > m$ and
$\gamma_2 >0, \gamma_2 >0$ so fine that if $(h_1,\ldots, h_n)  \in
\Gamma_R(z_1,\ldots, z_n;m_2,k,\gamma_2),$ then $f(h_1,\ldots, h_n) \in
\Gamma_L(z;m_1,k,\gamma_1).$

By the assumption for $k$ large enough there exists an $(h_{1k},\ldots,
h_{nk}) \in \Gamma_R(z_1,\ldots, z_n;m_2,k,\gamma_2).$ Thus, $x_k =
f(h_{1k},\ldots, h_{nk}) \in \Gamma_L(z;m_1,k,\gamma_1).$ Now recall the
matrices $y_k$ constructed in Lemma 6.2.  For $k$ large enough both $x_k$
and $y_k$ lie in $\Gamma_L(z;m_1,k,\gamma_1)$ and thus by the preceding
paragraph there exists a unitary $u$ satisfying $|uy_ku^* - x_k|_2 < t.$
It follows that if $\lambda_{ik}$ and $\mu_{ik}$ denote the respective
$i$th eigenvalues of $y_k$ and $x_k$ for $1 \leq i \leq k$ where the
eigenvalues are listed according from least to greatest and with respect
to multiplicity, then, $\sum_{i=1}^n |\lambda_{ik} - \mu_{ik}|^2 < t^2 k.$
With this in mind, $x_k$ is unitarily equivalent to the diagonal $k\times
k$

\[
\begin{bmatrix}
A^{\prime}_k & 0 \\
0 &  B^{\prime}_k \\
\end{bmatrix}
\]

\noindent where $A_k^{\prime}$ is a $p_k \times p_k$ diagonal matrix
($p_k$ defined in Lemma 6.2) and $B^{\prime}_k$ is a $(k-p_k) \times (k-
p_k)$ diagonal matrix and if $\lambda_{ik}$ is the $j^{th}$ eigenvalue of
the matrix above, then $\mu_{ik}$ is the $j^{th}$ eigenvalue of

\[ y_k =
\begin{bmatrix}
A_k & 0 \\
0 &  B_k \\
\end{bmatrix}
\]

\noindent We conclude that the $p_k \times p_k$ matrices $A_k$ and 
$A_k^{\prime}$ differ in $| \cdot|_2$-norm (on $M_{p_k}(\mathbb C)$) by 
no more than $t.$  

We're now going to compare the volumes of the $\epsilon$ neighborhoods
of the unitary orbits of $A_k$ and $A^{\prime}_k$ which we denote by
$\Theta_{\epsilon}(A_k)$ and $\Theta_{\epsilon}(A^{\prime}_k),$
respectively.  Denote again by $a_{1k}, \ldots, a_{p_kk}$ the
eigenvalues of $A_k$ ordered from least to greatest and according to
multiplicity and similarly denote by $a^{\prime}_{1k}, \ldots,
a^{\prime}_{p_kk}$ the eigenvalues of $A^{\prime}_k$ ordered in the 
same fashion.  We have $\sum_{i=1}^{p_k} |a_{ik} - a^{\prime}_{ik}|^2 < 
(\epsilon_0 \tau)^2 p_k.$ Now just as in the proof of Lemma 6.2 
$\text{vol}(\Theta_{\epsilon}(A^{\prime}_k))$ exceeds

\begin{eqnarray*} D_{p_k} \cdot \int_{\Omega^{\prime}_k} \Pi_{i <j}(t_i 
-t_j)^2 
\,dt_1\cdots dt_{p_k}. \end{eqnarray*}

\noindent where $\Omega^{\prime}_k$ is the intersection of

\[ [a^{\prime}_{1k} - \epsilon, a^{\prime}_{1k} + \epsilon] \times
\cdots \times [a^{\prime}_{p_k k} - \epsilon, a^{\prime}_{p_k k} +
\epsilon] \]

\noindent with $\mathbb R^{p_k}_{<}.$ Now we run the same argument in
Lemma 6.2.  Write $W_k^{\prime}$ for the set of 2-tuples (i,j) such that
$1 \leq i < j \leq p_k$ and $|a^{\prime}_{ik} - a^{\prime}_{jk}| <
\frac{\epsilon_0}{2}.$ Now consider $W_k^{\prime}-W_k.$ If $(i,j)$ is an
element of this set, then either $|a_{ik} - a^{\prime}_{ik}|$ or
$|a_{jk} - a_{jk}^{\prime}|$ exceeds $\frac{\epsilon_0}{4}.$  A moment's 
thought now shows that there exist at least $\#(W_k^{\prime} - W_k) / 
p_k$ indices $i, 1 \leq i \leq p_k,$ for which $|a_{ik} - 
a_{ik}^{\prime}| > \frac{\epsilon_0}{4}.$ Thus,

\[ (\epsilon_0 \tau)^2 \cdot p_k > \sum_{i=1}^{p_k} |a_{ik} - 
a_{ik}^{\prime}|^2 >
\frac{\#(W_k^{\prime} - W_k)}{p_k} \cdot \frac{\epsilon_0^2}{16}. \]

\noindent $\#(W^{\prime} - W_k) < 16 \tau^2 p_k^2.$  So $\#W^{\prime}_k 
< \#W_k + 16 \tau^2 p_k^2.$

Careful inspection of the chain of arguments in Lemma 6.2 now shows that
for all $0 < \epsilon < \epsilon_0,
\text{vol}(\Theta_{\epsilon}(A^{\prime}_k)) > L_k \epsilon^{(\beta + 16
\tau)k^2}.$ Again, Lemma 6.1 and the computations of the asymptotics of
$L_k$ show that for $\alpha_0 = 1 - 20 \tau - \sum_{j=1}^l c_j^2,$

\[ \liminf_{k\rightarrow \infty} k^{-2} \cdot H^{\alpha_0
k^2}_{\epsilon_0}(\Theta(x_k)) > \log \epsilon_0 - 17 \pi. \]

\noindent where $\Theta(x_k)$ is the unitary orbit of $x_k.$ Denote by
$\Theta(h_{1k}, \ldots, h_{nk})$ the set of all elements of the form
$(uh_{1k}u^*,\ldots, uh_{nk}u^*)$ where $u \in U_k$ and
observe that the map $f$ induces a function from $\Theta(h_{1k},\ldots,
h_{nk})$ onto $\Theta(x_k)$ and that this map has a Lipschitz constant 
(when both the domain and range are endowed with $|\cdot|_2$) $C$ where 
$C$ depends only upon $R$ and $f.$  It follows then that

\begin{eqnarray*} \mathbb H^{\alpha_0}_{\epsilon_0,R}(z_1,\ldots, z_n;m
,\gamma) & \geq & \liminf_{k \rightarrow \infty} k^{-2} \cdot
H_{\epsilon_0}^{\alpha_0 k^2} (\Theta(h_{1k},\ldots, h_{nk})) \\ & \geq &
\liminf_{k\rightarrow \infty}k^{-2}\cdot
H^{\alpha_0 k^2}_{\epsilon_0}(\Theta(x_k)) - \alpha_0 \log C \\ & > & 
\log\epsilon_0 - 17 \pi - \log C. \end{eqnarray*}

\noindent $m$ and $\gamma$ being arbitrary $- \infty < \mathbb
H^{\alpha_0}_{\epsilon_0, R}(z_1,\ldots,z_n) \leq \mathbb
H^{\alpha_0}_R(z_1,\ldots, z_n)  \leq \mathbb
H^{\alpha_0}(z_1,\ldots,z_n).$ Now $\alpha_0 = 1 -20 \tau -
\sum_{j=1}^{l}c_j^2$ and moreoever, $l \rightarrow s$ and $\alpha_0
\rightarrow 1 - \sum_{j=1}^s c_j^2 = \mathbb H(z)$ as $\tau \rightarrow
0.$ it follows that $\mathbb H^{\alpha}(z_1,\ldots, z_n) > - \infty$ for
all $\alpha < \mathbb H(z).$ Thus for such $z, \mathbb H(z_1,\ldots,
z_n) \geq \mathbb H(z).$

Finally suppose $z \in \{z_1, \ldots, z_n\}^{\prime \prime}.$  Find a 
sequence $\langle b_k \rangle_{k=1}^{\infty}$ such that $b_k \rightarrow 
z$ strongly.  For each $k$ the preceding argument shows that $\mathbb 
H(z_1,\ldots, z_n) \geq \mathbb H (b_k)$ whence by Corollary 7.2 $\mathbb 
H(z_1,\ldots, z_n) \geq \liminf_{k \rightarrow \infty} \mathbb H(b_k) \geq 
\mathbb H(z).$
\end{proof}

\begin{corollary} If $\{z_1,\ldots, z_n\}$ has finite dimensional 
approximants and $\mathbb H(z_1,\ldots, z_n) <1,$ then $\{z_1,\ldots, 
z_n \}^{\prime \prime}$ has a minimal projection. 
\end{corollary}

\begin{proof} Suppose $\{z_1, \ldots, z_n \}^{\prime \prime}$ is diffuse,
i.e., has no minimal projections.  Find a maximal abelian subalgebra $N$
of $\{z_1, \ldots, z_n \}^{\prime \prime}$ and a single selfadjoint
generator $z$ for $N.$ $z$ has no eigenvalues by maximality of $N$ and
thus by Lemma 6.4 $\mathbb H(z_1,\ldots, z_n) \geq \mathbb H(z) =1.$
\end{proof}

\section{Additivity Properties of $\mathbb H$}

In this section we prove additive formulae for $\mathbb H$ in the
presence of freeness.  

\begin{theorem} If $\{z_1,\ldots,z_n\}$ is set of freely independent,
selfadjoint elements of $M,$ then

\[ \mathbb H(z_1, \ldots,z_n) = \sum_{i=1}^n \mathbb H(z_i).\]
\end{theorem}

\begin{proof} Observe by Corollary 3.5 and Lemma 6.2 that  

\[ \mathbb H(z_1,\ldots,z_n) \leq \delta_0(z_1,\ldots,z_n) \leq
\sum_{i=1}^n \delta_0(z_i)  = \sum_{i=1}^n\mathbb H(z_i) .\]

\noindent Thus it remains to show that $\mathbb H(z_1,\ldots,z_n) \geq
\sum_{i=1}^n\mathbb H(z_i).$  

	For $1 \leq i \leq n$ define $\alpha_i = \mathbb H(z_i).$ Set
$\alpha = \alpha_1 + \cdots + \alpha_n.$ Suppose $m
\in \mathbb N,$ $\tau, \gamma >0,$ and $R > \max \{ \|z_i\|\}_{1 \leq i 
\leq n}.$  By Corollary 2.14 of [11] there exists an
$N\in \mathbb N$ such that if $k \geq N$ and $\sigma$ is a Radon
probability measure on $((M_{k}^{sa}(\mathbb C))_{R})^{n}$ invariant
under the $ (U_k)^{(n-1)} $-action given by $(\xi_1,\ldots,\xi_n)  
\mapsto (\xi_1, u_1\xi_2u^*_1, \ldots, u_{n-1}\xi_nu^*_{n-1})$ where $
(u_1,\ldots,u_{n-1}) \in (U_k)^{(n-1)},$ then $\sigma(\omega_k) > \frac
{1}{2}$ where 

\[ \omega_k = \{(\xi_1,\ldots,\xi_n)\in((M_k^{sa}(\mathbb 
C))_{R+1})^{n} : \langle \{\xi_i\} \rangle_{i=1}^n \hspace{.05in} 
\text{are} \left (m, \frac{\gamma}{4^m} \right)
\text{ - free} \}.
\]

\noindent The preceding section provided for each $i$ a sequence $\langle
y_{ik} \rangle_{k=1}^{\infty}$ such that for any $m^{\prime} \in \mathbb
N$ and $\gamma^{\prime} >0$ $y_{ik} \in
\Gamma_R(z_i;m^{\prime},k,\gamma^{\prime})$ for sufficiently large $k.$
Write $\Theta (y_{ik})$ for the unitary orbit of $y_{ik}$ and
$g_{ik}$ for the topological dimension of this orbit.  The proof of
Lemma 6.2 yields a $1 > \epsilon_0 >0$ such that for each $i$ and $k$
sufficiently large there exist constants $L_{ik}$ and $b_{ik}, 
p_{ik} \in \mathbb
N$ such that for $\epsilon_0 > \epsilon >0$

\[ P_{\epsilon}(\Theta(y_{ik})) \geq L^{(i)}_k
\cdot \frac{\epsilon^{(b_{ik} - p_{ik}^2)} 
\Gamma(\frac{p_{ik}^2}{2} +1 )}{4^{k^2} \sqrt{\pi k}^{p_{ik}^2} }
\]

\noindent We may arrange it so that if $\beta_k = 
\sum_{i=1}^n p_{ik}^2 - b_{ik},$ then for sufficiently large $k$ 
$\beta_k > (\alpha- \tau) k^2.$
 
For each $k \in \mathbb N$ denote by $\mu_k$ the probability measure
on $((M^{sa}_k(\mathbb C))_{R+1})^n$ obtained by restricting $\sum_{i=1}^n
g_{ik}$-Hausdorff measure (with respect to the $| \cdot|_2$ norm) to
the smooth $\sum_{i=1}^n g_{ik}$-dimensional manifold $T_k =
\Theta(y_{1k}) \times \cdots \times \Theta(y_{nk})$ and normalizing
appropriately.  $\mu_k$ is a Radon probability measure invariant under the
$(U_k)^{n-1}$-action in the sense described above (such an action is
isometric and thus does not alter Hausdorff measure).  $\mu_k(\omega_k) >
\frac{1}{2}.$ Define $F_k = \omega_k \bigcap T_k.$ It is clear that
$\mu_k(F_k)= \mu_k(\omega_k) > \frac{1}{2}$ and $F_k \subset
\Gamma_{R+1}(z_1,\ldots,z_n;m,k,\gamma).$ It remains to make lower bounds
on the Hausdorff quantities of $F_k.$
  
  $T_k$ is a locally isometric smooth manifold of 
dimension $\sum_{i=1}^n g_{ik}.$  From the preceding paragraph it 
follows that for all $0 < \epsilon < \epsilon_0 $

\[P_{\epsilon}(T_k) \geq \Pi_{i=1}^n P_{\epsilon} (\Theta(y_{ik})) \geq
\Pi_{i=1}^n L_{ik} \cdot \left (\frac{\epsilon^{(b_{ik} -
p_{ik}^2)} \Gamma(\frac{p_{ik}^2}{2} +1 )}{4^{k^2}
\sqrt{\pi k}^{p_{ik}^2} } \right). \] 

\noindent By Lemma 6.1

\[ H^{\beta_k}_{\epsilon_0} (F_k) > \Pi_{i=1}^n \cdot \frac{ 
L_{ik}
\Gamma(\frac{p_{ik}^2}{2} +1 )}{4^{k^2} \sqrt{\pi
k}^{p_{ik}^2} } \cdot \frac{1}{2}. \]

\noindent For any $m \in \mathbb N, \gamma >0$ $\mathbb H^{\alpha -
\tau}_{\epsilon_0, R}(z_1,\ldots,z_n;m,\gamma)$ dominates 

\begin{eqnarray*}
\liminf_{k \rightarrow \infty} k^{-2} \log
(H^{\beta_k}_{\epsilon_0}(F_k)) & \geq & 
\liminf_{k \rightarrow 
\infty} 
k^{-2} \log \left ( \Pi_{i=1}^n 
\left (\frac{L_{ik} \Gamma(\frac{p_{ik}^2}{2} +1 
)}{4^{k^2} \sqrt{\pi
k}^{p_{ik}^2} } \right) \cdot \frac{1}{2} \right)
\\ & \geq & 
\liminf_{k \rightarrow \infty} \sum_{i=1}^n \left [k^{-2} \left (\log 
L_{ik}
\Gamma \left( \frac {p_{ik}^2}{2} +1\right )\right) - 
\frac{p_{ik}^2}{2k^2} \log k  
- 16 \pi \right ]  \\
& \geq & 
\sum_{i=1}^n \liminf_{k\rightarrow \infty} \left[ \left(k^{-2} \log       
\left ( L_{ik}
\Gamma \left (\frac {p_{ik}^2}{2} +1\right)\right) - 
\frac{p_{ik}^2}{2k^2} \log k \right)        
- 16 \pi \right]\\ & > & n \cdot \log \epsilon_0 - 17n \pi.
\end{eqnarray*} 

\noindent $\mathbb H^{\alpha - \tau}(z_1,\ldots,z_n)  \geq
\mathbb H^{\alpha -\tau}_{\epsilon_0,R}(z_1,\ldots,z_n) > -\infty.$
$\tau>0$ being arbitrarily small $\mathbb H(z_1,\ldots,z_n) \geq \alpha
= \mathbb H(z_1) + \cdots + \mathbb H(z_n).$ \end{proof}

We now turn to the situation where we have free products of finite
dimensional algebras.  We obtain a slightly stronger result.  The
arguments proceed as in Theorem 7.1 but the issues are a bit more
delicate.  First we rephrase Lemma 5.4 in terms of $\epsilon$ packings.

\begin{lemma} Suppose $\{z_1,\ldots,z_n\}$ generates a finite dimensional
unital subalgebra $A$ of $M,$ $\alpha = \delta_0(A),$ and $R > \max \{
\|z_i\| \}_{1 \leq i \leq n}.$ There exists constants $K >0 $ and
$\epsilon_0>0$ such that for any given $m\in \mathbb N$ and $\gamma >0$
there exists an integer $N$ such that for $k \geq N$ there is a locally
isometric smooth manifold $T_k \subset
\Gamma_R(z_1,\ldots,z_n;m,k,\gamma)$ of dimension $g_k \geq \alpha k^2$
and for any $0 < \epsilon < \epsilon_0$

\[ P_{\epsilon}(T_k) \geq \left( \frac 
{K}{\epsilon}\right)^{\alpha k^2}\] \end{lemma}

\begin{proof} By Lemma 5.4 there exist $\lambda, r, \zeta, D, L> 0, N \in
\mathbb N$ such that for any $\epsilon < \lambda(DL)^{-1}, m \in \mathbb
N, \gamma >0$ and any $k \geq N$ there exists a locally isometric
smooth manifold $T_k$ (obtained by smearing the images of the $z_i$ under
a representation with $U_k$) with $T_k \subset 
\Gamma_R(z_1,\ldots,z_n;m,k,\gamma)$ and

\[ H^{\alpha k^2}_{\epsilon}(T_k) \geq \left( \frac{\lambda r \zeta}{DL} 
\right)^{\alpha k^2}.\]

\noindent Thus 

\[ P_{\frac{\epsilon}{2}}(T_k) \cdot \epsilon^{\alpha k^2} \geq \left( \frac{\lambda r \zeta}{DL} \right)^{\alpha
k^2}\]

\noindent whence,

\[ P_{\frac{\epsilon}{2}}(T_k) \geq \left( \frac{\lambda r 
\zeta}{DL\epsilon} \right)^{\alpha k^2} 
.\]
\noindent Set $K = \frac{\lambda r \zeta}{2DL}.$ \end{proof} 

Suppose $Z_1, \ldots Z_n$ are finite ordered sets of selfadjoint elements
in $M.$ We write $Z_1 \cup \cdots \cup Z_n$ for the ordered set obtained
by listing the elements of $Z_1$ in order, then $Z_2,$ etc.  It is in this
way that we interpret $\Gamma_R (Z_1 \cup \ldots \cup Z_n;m,k,\gamma)$ and 
all the asymptotic dimensions and measurements associated to
$Z_1 \cup \ldots \cup Z_n.$

\begin{lemma} If $Z_1,\ldots, Z_n$ are ordered sets of selfadjoint 
elements in $M$ with $\alpha_i = \delta_0(Z_i)$ and $\alpha = \alpha_1 + 
\cdots +\alpha_n,$ then  

\[ \mathbb P^{\alpha}(Z_1 \cup \ldots \cup Z_n) \leq \mathbb 
P^{\alpha_1}(Z_1) + 
\cdots + \mathbb P^{\alpha_n}(Z_n) + \alpha \log(4 \sqrt{n}).\]

\end{lemma}

\begin{proof} For any $R, \gamma >0$ and $m,k\in \mathbb N$

\[ \Gamma_R(Z_1 \cup \ldots \cup Z_n;m,k,\gamma) \subset 
\Gamma_R(Z_1;m,k,\gamma) 
\times \cdots \times \Gamma_R(Z_n;m,k,\gamma) \]

\noindent The proof now follows from going through the definitions 
and using subadditivity of $P_{\epsilon}$ on products.  We have

\begin{eqnarray*} P_{4 \epsilon 
\sqrt{n} }(\Gamma_R(Z_1 \cup \ldots \cup Z_n;m,k,\gamma) & \leq 
& P_{4 \epsilon \sqrt{n}}(\Gamma_R(Z_1;m,k,\gamma) \times \cdots \times 
\Gamma_R(Z_n;m,k,\gamma)) \\ & \leq & \Pi_{i=1}^n 
P_{\epsilon}(\Gamma_R(Z_i;m,k,\gamma)).
\end{eqnarray*}
 
\noindent Thus $\mathbb P_{4 \epsilon \sqrt{n} 
,R}(Z_1 \cup \ldots \cup Z_N;m,\gamma) \leq 
\sum_{i=1}^n \mathbb P_{\epsilon,R}(Z_i;m,\gamma) \Rightarrow 
\mathbb P_{4 \epsilon \sqrt{n} }(Z_1,\ldots,Z_n) \leq \sum_{i=1}^n \mathbb 
P_{\epsilon}(Z_i).$  Consequently,

\begin{eqnarray*} \mathbb P^{\alpha}(Z_1 \cup \ldots \cup Z_n) & = & 
\limsup_{\epsilon 
\rightarrow 0} 
\left (\mathbb 
P_{4\epsilon \sqrt{n}}(Z_1,\ldots,Z_n) 
+ \alpha \log (8\epsilon \sqrt{n})  \right) \\
& \leq & \sum_{i=1}^n \limsup_{\epsilon \rightarrow 0} \left (\mathbb 
P_{\epsilon}(Z_i) + \alpha_i \log (8 \epsilon \sqrt{n}) \right)\\ 
& = & \mathbb P^{\alpha_1}(Z_1) + \cdots + \mathbb P^{\alpha_n}(Z_n) + 
\alpha \log (4\sqrt{n}). 
\end{eqnarray*}

\end{proof}

\begin{theorem} Suppose $\{ Z_1,\ldots, Z_n \}$ is a freely independent
family of ordered sets of selfadjoint elements in $M$ such that each $Z_i
= \{ z_{i1},\ldots, z_{ip_i}\}$ generates a finite dimensional
unital subalgebra $A_i.$ If for each $i$ $\alpha_i = \mathbb H(A_i)$ and
$\alpha = \alpha_1 + \cdots + \alpha_n,$ then $Z_1 \cup \ldots \cup Z_n$ 
is an $\alpha$-set and

\[ \mathbb H(Z_1 \cup \ldots \cup Z_n) = \sum_{i =1}^n \mathbb H(Z_i) = 
\sum_{i=1}^n \delta_0(Z_i) .\]
\end{theorem}
\begin{proof}  It suffices to show that $\{Z_1,\ldots,Z_n\}$ is an 
$\alpha$-set for the second assertion is an immediate consequence of this.

	First we show that $\mathbb H^{\alpha}(Z_1 \cup \ldots \cup Z_n) >
-\infty.$ We will write elements of $((M_k^{sa}(\mathbb C))^{p_1 + \ldots
+ p_n}$ as $(X_1,\ldots,X_n)$ where $X_i \in (M_{k}^{sa}(\mathbb
C))^{p_i}.$ Moreover, for any $u \in U_k$ $uX_iu^*$ denotes the element of
$(M^{sa}_k(\mathbb C))^{p_i}$ obtained by conjugating each entry of $X_i$
by $u.$ Suppose $m \in \mathbb N,$ $\gamma >0,$ and $R$ exceeds the
operator norm of any element in one of the $Z_i.$ Again Corollary 2.14 of
[11] yields an $N\in \mathbb N$ such that if $k \geq N$ and $\sigma$ is a
Radon probability measure on $((M_k^{sa}(\mathbb C))_{R+1})^{p_1 + \ldots
+ p_n}$ invariant under the $ (U_k)^{(n-1)} $-action given by

\[ (X_1,\ldots, X_n) \mapsto (X_1, u_1X_2u_1^*,\ldots, 
u_{n-1}X_nu_{n-1}^*)
\]

\noindent where $(u_1,\ldots,u_{n-1}) \in (U_k)^{(n-1)},$ then 
$\sigma(\omega_k) > \frac {1}{2}$ where

\[ \omega_k = \{(X_1,\ldots, X_n)\in((M_k^{sa}(\mathbb
C))_{R+1})^{p_1 + \cdots + p_n} : \langle X_i \rangle_{i=1}^n
\hspace{.05in}
\text{are} \left (m, \frac{\gamma}{4^m} \right)
\text{ - free} \}.
\]

By Lemma 7.2 there exists for each $i$ constants $K_i>0$ and $\epsilon_i
>0$ such that for a given $m \in \mathbb N$ and $\gamma >0$ there exists
$N_i \in \mathbb N$ such that for $k \geq N_i$ there exists a locally
isometric smooth manifold $T_{ik}$ of dimension $g_{ik} \geq \alpha_i
k^2$ with $T_{ik} \subset \Gamma_R\left(Z_i;m,k,
\frac{\gamma}{(8(R+1))^m}\right)$ and for any $0 < \epsilon < \epsilon_i$

\[ P_{\epsilon}(T_k) \geq \left(\frac {K_i}{\epsilon} \right)^{\alpha_i 
k^2}.\]

\noindent For any $k > N_1 + \cdots + N_n$ set $\Omega_k = T_k^1
\times \cdots \times T_k^n$ and denote by $\mu_k$ the probability measure
on $((M^{sa}_k(\mathbb C))_{R+1})^{p_1+ \cdots + p_n}$ obtained by
restricting $\sum_{i=1}^n g_{ik}$-Hausdorff measure (with respect to
the $| \cdot |_2$ norm) to $\Omega_k$ and normalizing appropriately.  As
in Theorem 7.1 $\mu_k$ is a Radon probability measure invariant under the
$(U_k)^{n-1}$ -action described above, whence $\mu_k(\omega_k) >
\frac{1}{2}.$ Define $F_k = \omega_k \bigcap \Omega_k.$ $\mu_k(F_k) =
\mu_k(\omega_k) > \frac{1}{2}$ and $F_k \subset
\Gamma_{R+1}(Z_1 \cup \ldots \cup Z_n;m,k,\gamma).$ 

For each $k$ sufficiently large $\Omega_k$ is a locally
isometric smooth manifold of dimension $g_k = \sum_{i=1}^n g_{ik}$.  
Moreover setting $K = \min \{K_i \}_{1\leq i \leq n}$ for all $0 < 
\epsilon < \min_{1 \leq i \leq n} \{\epsilon_i\}$

\[ P_{\epsilon}(\Omega_k) \geq \Pi_{i=1}^n 
P_{\epsilon} (T_k^i) \geq \left(\frac 
{K}{\epsilon}\right)^{\alpha k^2}
\]
\noindent By Lemma 6.1 for $\epsilon_0 = \min_{1 \leq i \leq n} 
\{\epsilon_i\}$

\[ H^{\alpha k^2}_{\epsilon_0}(F_k) > 
\left(\frac{K}{2}\right)^{\alpha k^2} \cdot \frac {1}{2}.\] 

\noindent Thus,

\[ k^{-2} \cdot \log H^{\alpha k^2}_{\epsilon_0}(F_k) > \alpha \log K - 
\log 4.\]

Given $m \in \mathbb N$ and $\gamma >0$ there exists for each $k$ 
large enough a set $F_k \subset 
\Gamma_R(Z_1 \cup \ldots \cup Z_n;m,k,\gamma)$ satisfying 
the outer Hausdorff measure lower bound above.  Consequently, $\mathbb 
H^{\alpha}_{\epsilon_0,R}(Z_1 \cup \ldots \cup Z_n;m,\gamma)$ dominates

\[ \alpha \log K -4 > -\infty
\]

\noindent This lower bound is independent of $m$ and $\gamma$ so it 
follows that the above expression is a lower bound for $\mathbb 
H^{\alpha}_{\epsilon_0}(Z_1 \cup \ldots \cup Z_n).$  It follows that 
$\mathbb H^{\alpha}(Z_1 \cup \ldots \cup Z_n) > -\infty$ as promised.  

It remains to show that $\mathbb H^{\alpha}(Z_1 \cup \ldots \cup Z_n) < 
\infty.$ We have by Lemma 7.3 and Lemma 5.1 that

\[ \mathbb H^{\alpha}(Z_1 \cup \ldots \cup Z_n) - \alpha \log (8 \sqrt{n})  
\leq \mathbb P^{\alpha}(Z_1 \cup \ldots \cup Z_n) - \alpha \log(4
\sqrt{n})  \leq \sum_{i=1}^n \mathbb P^{\alpha_i}(Z_i) < \infty.\]
\end{proof}

\noindent{\it Acknowledgements.} Part of this research
was conducted during the 2002 PIMS conference on asymptotic geometry in
Vancouver.  I thank Stanislaw Szarek and Gilles Pisier for their
hospitality.  The former made a helpful remark on the
relation between Hausdorff dimension and Minkowski content on smooth
manifolds.  Another part of this research took place at UCLA and I thank
Dimitri Shlyakhtenko and the UCLA mathematics department for their
hospitality.

\end{document}